%% file: tex-me.tex
\documentclass[11pt]{amsart}
\newcommand\version{public}

\input{initializing}
\usepackage{youngtab}
\usepackage{amsmath,amsthm,xspace,xcolor}
\usepackage{cite,pgf,pgfarrows,pgfnodes,pgfautomata,pgfheaps}
\usepackage{tikz}
\usepackage{placeins}
\usepackage{hyperref}

\input{definitions}

\begin{document}
\input{article-title}

\input{hyp10max}

\input{tex-me.bbl}
\end{document}

%% file: initializing.tex
\usepackage{fullpage,setspace}
\usepackage{enumerate}   
\usepackage{url}
\usepackage{enumerate}
\usepackage{amsmath,amssymb,amscd,mathrsfs,latexsym}
\usepackage{mathdots}

\usepackage{ifthen}     

\newcommand\pubpri[2]{%
\ifthenelse{\equal{\version}{public}}%
{{#1}}%
{\marginpar{\scshape\small Pubpri Alert}{#2}}}
\newcommand\pubprinoalert[2]{%
\ifthenelse{\equal{\version}{public}}%
{{#1}}%
{#2}}

\newcommand\ignore[1]{}

\usepackage{color}
\providecommand\wantcolor{yes}   %
\definecolor{backgroundyellow}{cmyk}{.2,.1,.8,.2}
\definecolor{backgroundblue}{rgb}{0,0,1}
\definecolor{backgroundred}{rgb}{1,0,0}
\definecolor{backgroundmagenta}{cmyk}{0,1,0,0}
\ifthenelse{\equal{\wantcolor}{no}}{\renewcommand\color[1]{}}{}

\newcommand\mysubsubsection[1]{%
		\subsubsection{\sffamily\upshape\mdseries #1}}

\newcommand\mysss{\mysubsubsection}

\usepackage{chngcntr}
\counterwithin{figure}{section}

\providecommand{\theoremnumbering}{section}

\ifthenelse{\equal{\theoremnumbering}{section}}{\newtheorem{annotation}{Annotation}[section]}%
{\ifthenelse{\equal{\theoremnumbering}{subsection}}{\newtheorem{annotation}{Annotation}[subsection]}{\newtheorem{annotation}{Annotation}}}
\newtheorem{theorem}[annotation]{
		Theorem}
\newtheorem{lemma}[annotation]{
		Lemma}
\newtheorem{definition}[annotation]{
		Definition}
\newtheorem{corollary}[annotation]{
		Corollary}
\newtheorem{proposition}[annotation]{
		Proposition}

\newtheorem{example}[annotation]{
		Example}
\newcommand\bexample{\begin{example}\begin{rm}}
\newcommand\eexample{\end{rm}\hfill$\Box$\end{example}}
\newtheorem{examplenobox}[annotation]{
		Example}
\newcommand\bexamplenobox{\begin{examplenobox}\begin{rm}}
\newcommand\eexamplenobox{\end{rm}\end{examplenobox}}
\newtheorem{exercise}[annotation]{
		Exercise}
\newcommand\bexercise{\begin{exercise}\begin{rm}}
\newcommand\eexercise{\end{rm}\end{exercise}}
\newtheorem{notation}[annotation]{
		Notation}
\newcommand\bnotation{\begin{notation}\begin{rm}}
\newcommand\enotation{\end{rm}\end{notation}}

\newtheorem{remark}[annotation]{
		Remark}
\newcommand\bremark{\begin{remark}
\begin{upshape}}
\newcommand\eremark{\end{upshape}
\end{remark}}
\newenvironment{remark*}{%
\par\noindent{\scshape 
  Remark: }\begin{rm}}{\hfill\end{rm}\newline} 
\newcommand\bremarkstar{\begin{remark*}}
\newcommand\eremarkstar{\end{remark*}}
\newcommand\bdefn{\begin{definition}
\begin{upshape}}
\newcommand\edefn{\end{upshape}
\end{definition}}
\newtheorem{caveat}[annotation]{
		Caveat}
\newcommand\bcaveat{\begin{caveat}
\begin{upshape}}
\newcommand\ecaveat{\end{upshape}
\end{caveat}}
\newenvironment{caveatstar}{
\par\noindent{\scshape\bfseries
  Caveat: }\begin{rm}}{\end{rm}\newline} 
\newcommand\bcaveatstar{\begin{caveatstar}}
\newcommand\ecaveatstar{\end{caveatstar}}
\newenvironment{myproof}{%
\par\noindent{\scshape 
  Proof: }\begin{rm}}{\hfill$\Box$\end{rm}\newline} 
\newcommand\bmyproof{\begin{myproof}}
\newcommand\emyproof{\end{myproof}}
\newenvironment{myproofnobox}{%
\par\noindent{\scshape Proof: }\begin{rm}}{\end{rm}\hfill\newline}
\newcommand\bmyproofnobox{\begin{myproofnobox}}
\newcommand\emyproofnobox{\end{myproofnobox}}
\newenvironment{solution}{%
\par\noindent{\scshape Solution: }\begin{rm}}{\hfill$\Box$\end{rm}\newline}
\newenvironment{solutionnobox}{%
\par\noindent{\scshape Solution: }\begin{rm}}{\end{rm}}
\newcommand\bsolution{\begin{solution}\begin{rm}}
\newcommand\esolution{\end{rm}\end{solution}}
\newcommand\bsolutionnobox{\begin{solutionnobox}\begin{rm}}
\newcommand\esolutionnobox{\end{rm}\end{solutionnobox}}
\newcommand\bthm{\begin{theorem}}
\newcommand\ethm{\end{theorem}}
\newcommand\bcor{\begin{corollary}}
\newcommand\ecor{\end{corollary}}
\newcommand\blemma{\begin{lemma}}
\newcommand\elemma{\end{lemma}}
\newcommand\bprop{\begin{proposition}}
\newcommand\eprop{\end{proposition}}
\newcommand\beqn{\begin{equation}}
\newcommand\eeqn{\end{equation}}
\newcommand\beqnstar{\begin{equation*}}
\newcommand\eeqnstar{\end{equation*}}
\numberwithin{equation}{section}

\newcommand\mtitle[1]%
{\marginpar{\sffamily\raggedright\upshape\small #1}}

\providecommand\finalized{yes}
\ifthenelse{\equal{\finalized}{yes}}{}{\marginparwidth=2cm}
\ifthenelse{\equal{\finalized}{yes}}%
{\newcommand\mylabel[1]{\label{#1}}}%
{\newcommand\mylabel[1]{\label{#1}\marginpar{[{\ttfamily\upshape\tiny #1}]}}}
\ifthenelse{\equal{\finalized}{yes}}%
{\newcommand\checked[1]{}}%
{\newcommand\checked[1]{\marginpar{[{\ttfamily\upshape\tiny CHECKED: #1}]}}}
\ifthenelse{\equal{\finalized}{yes}}%
{\newcommand\spellchecked[1]{}}%
{\newcommand\spellchecked[1]{\marginpar{[{\ttfamily\upshape\tiny SPELLCHECKED: #1}]}}}

\providecommand\version{public}   
\ifthenelse{\equal{\version}{public}}%
{\newcommand\mcomment[1]{}}%
{\newcommand\mcomment[1]{\marginpar{{\raggedright\sffamily\upshape\small
\begin{spacing}{0.75} #1\end{spacing}}}}}
\ifthenelse{\equal{\version}{public}}%
{\newcommand\fcomment[1]{}}%
{\newcommand\fcomment[1]{\footnote{#1}}}
\ifthenelse{\equal{\version}{public}}%
{\newcommand\comment[1]{}}%
{\newcommand\comment[1]{{\small #1}}}

\usepackage{graphicx}
\graphicspath{ {diagram/} }


%% file: definitions.tex
\newcommand\algA{\lieg(A)}
\newcommand\algB{\lieg(B)}
\newcommand{\algprimeA}[1][A]{\liegprime(#1)}
\newcommand\algprimeB{\algprimeA[B]}
\newcommand\deralgB{\liegprime(B)}
\newcommand\lieg{\mathfrak{g}}
\newcommand\liegprime{\mathfrak{g}'}
\newcommand\lieh{\mathfrak{h}}
\newcommand\liek{\mathfrak{k}}
\DeclareMathOperator{\rank}{\textup{rank}}

\newcommand\pisys{$\pi$-system\xspace}
\newcommand\pisystem{\pisys}
\newcommand\pisystems{$\pi$-systems\xspace}
\newcommand{\qsig}[1][\Sigma]{q_{{}_{#1}}}
\newcommand{\isig}[1][\Sigma]{i_{{}_{#1}}}
\newcommand{\Xsh}[1][X]{X_{\mathrm{short}}}
\newcommand{\Xlong}[1][X]{X_{\mathrm{long}}}

\DeclareMathOperator{\ad}{ad}

\newcommand{\roots}[1][]{\Delta_{#1}}
\newcommand{\reroots}[1][]{\Delta^{re}_{#1}}
\newcommand{\imroots}[1][]{\Delta^{im}_{#1}}

\newcommand{\Ext}{\mathrm{Ext}}

\newcommand{\Hyp}{\mathrm{Hyp}}

\newcommand{\finpart}[1]{{#1}^\circ}
\newcommand{\affpart}[1]{\widehat{\finpart{#1}}}
\newcommand{\form}[2]{\left( #1 \mid #2 \right)}
\newcommand{\aform}[2]{\langle\, #1, #2 \,\rangle}
\newcommand{\ovex}[1]{{#1}^{++}} 
\newcommand{\hasl}{\ovex{A}_1}
\newcommand{\afsl}{A^{(1)}_1}

\newcommand{\integers}{\mathbb{Z}}
\newcommand{\real}{\mathbb{R}}
\newcommand{\complex}{\mathbb{C}}
\newcommand{\euc}[1][]{\mathbb{E}_{#1}}

\newcommand{\cim}{C^{\mathrm{im}}}
\newcommand{\cimcl}{\overline{\cim}}
\newcommand{\cre}{C^{\mathrm{re}}}

\newcommand{\norm}[1]{\Vert #1 \Vert}
\newcommand{\be}{\begin{enumerate}}
\newcommand{\ee}{\end{enumerate}}
\DeclareMathOperator{\supp}{supp}
\DeclareMathOperator{\mult}{\mathfrak{m}}


%% file: article-title.tex
\author{Lisa Carbone}
\address{Department of Mathematics, Rutgers, Piscataway, NJ 08854-8019, USA}
\email{carbonel@math.rutgers.edu}

\author{K. N. Raghavan}
\address{The Institute of Mathematical Sciences, HBNI, Chennai 600113, India}
\email{knr@imsc.res.in}

\author{Biswajit Ransingh}
\address{Harish-Chandra Research Institute, HBNI, Allahabad 211019, India}
\email{biswajitransingh@hri.res.in}

\author{Krishanu Roy}
\address{The Institute of Mathematical Sciences, HBNI, Chennai 600113, India}
\email{krishanur@imsc.res.in}

\author{Sankaran Viswanath}
\address{The Institute of Mathematical Sciences, HBNI, Chennai 600113, India}
\email{svis@imsc.res.in}

\thanks{The first author's research is partially supported by the Simons Foundation, Mathematics and Physical Sciences-Collaboration Grants for Mathematicians, Award Number 422182. KNR and SV acknowledge support from DAE under a XII plan project. KR acknowledges SICI for an SRSF grant.}

\title{${\boldsymbol\pi}$-systems of symmetrizable Kac-Moody algebras}

\thanks{}

\begin{abstract}
As part of his classification of regular semisimple subalgebras of semisimple Lie algebras, Dynkin introduced the notion of a $\pi$-system. This is a subset of the set of roots such that pairwise differences of its elements are not roots. Such systems arise as simple systems of regular semisimple subalgebras.
Morita and Naito generalized this notion to all symmetrizable Kac-Moody algebras. In this work, we systematically develop the theory of $\pi$-systems of symmetrizable Kac-Moody algebras and establish their fundamental properties. For several Kac-Moody algebras with physical significance, we study the orbits of the Weyl group on $\pi$-systems, and completely determine the number of orbits.  In particular, we show that there is a unique $\pi$-system of type $\ovex{A}_1$ (the Feingold-Frenkel rank 3 hyperbolic algebra) in $E_{10}$ (the rank 10 hyperbolic algebra) up to Weyl group action and negation. 
\end{abstract}
\keywords{\pisystem, regular subalgebra, Weyl group orbits, symmetrizable Kac-Moody algebras}
\subjclass[2010]{17B22 (17B67)}
\maketitle

%% file: hyp10max.tex
\section{Introduction}

\subsection{}
Let $\lieg$ denote a finite dimensional semisimple Lie algebra over $\complex$. A semisimple Lie subalgebra of $\lieg$ is said to be regular if
it is $\ad \lieh$-invariant for some Cartan subalgebra $\lieh$ of $\lieg$.
The classical works of Borel-de Siebenthal \cite{bds} and Dynkin \cite{dynkin} contain a  complete classification of the possible Cartan-Dynkin types of regular semisimple subalgebras of each simple Lie algebra. 

Now let $\lieh$ be a fixed Cartan subalgebra and $\roots$ the corresponding set of roots. A \pisystem $\Sigma$ is a subset of $\roots$ satisfying the property that $\alpha - \beta \not\in \roots$ for all pairs of roots $\alpha \neq \beta$. Dynkin showed that linearly independent \pisystems arise precisely as simple systems of regular semisimple subalgebras of $\lieg$.
The Weyl group $W(\lieg) \subset \mathrm{GL}(\lieh^*)$  acts naturally on \pisystems, and the number of orbits for each type of \pisystem was tabulated by Dynkin in \cite[Tables 9-11]{dynkin}.


More generally, for a symmetrizable Kac-Moody algebra  $\lieg$, a \pisystem $\Sigma$ is a subset of its {\em real roots} such that pairwise differences of elements of $\Sigma$ are not roots of $\lieg$ \cite{morita, naito}. \pisystems of affine Kac-Moody algebras were classified in   \cite{naito, roy-venkatesh, FRT}. Unsurprisingly, much less is known beyond the affine case. However, many \pisystems of the hyperbolic Lie algebra $E_{10}$ have been constructed and studied for their significance in physics  \cite{feingold-nicolai, HLPS, Henneaux2008}. Much of the current paper stems from an attempt at understanding the $W(E_{10})$-orbits for the simplest \pisystem of $E_{10}$, namely that of type $\ovex{A}_1$ (the rank 3 hyperbolic Feingold-Frenkel algebra).

\subsection{} We now describe our principal results. Let $W$ denote the Weyl group of the symmetrizable Kac-Moody algebra $\lieg$. We first establish (Theorem~\ref{thm:pos-neg})  that every indecomposable linearly independent \pisystem $\Sigma$ of $\lieg$ is $W$-conjugate to one in which all roots have the same sign.
When $|\Sigma|=2$, this is a result of Naito \cite[Prop. 2.1]{naito}.

Consider the action of $W \times \integers_2$ on \pisystems of $\lieg$,  where  $\integers_2$ acts by $\Sigma \mapsto \pm \Sigma$.
If $\lieg$ has a  {\em symmetric} generalized Cartan matrix (GCM), we show that it admits a \pisystem of type $\ovex{A}_1$ if and only if its Dynkin diagram has a subdiagram of overextended type, i.e., one obtained by overextending some finite type Dynkin diagram (figure~\ref{fig:exts}). We establish a somewhat remarkable bijective correspondence between such subdiagrams of $\lieg$ and $W \times \integers_2$-orbits of \pisystems of type $\ovex{A}_1$ (Theorem~\ref{thm:ha11-pi-systems}).
It follows that that there is a unique \pisystem (up to $W \times \integers_2$-equivalence) of type $\ovex{A}_1$ in $E_{10}$. More generally, this uniqueness result holds when $E_{10}$ is replaced by any of the simply-laced overextended hyperbolic diagrams $\ovex{A}_n, \ovex{D}_n, \ovex{E}_n$ (corollary \ref{cor:hyp-a11}).

Strengthening these arguments, we prove a general counting theorem (Theorem~\ref{thm:mainthm-gen}) for $W \times \integers_2$-orbits of \pisystems of overextended type in $\lieg$. For a fixed overextended diagram $K$, the problem of counting the number of inequivalent \pisystems of type $K$ in $\lieg$ can be reduced to one in which $K$ is replaced by its underlying finite type diagram $\finpart{K}$ and $\lieg$ is replaced by certain finite-dimensional simple Lie subalgebras of $\lieg$ (equation~\eqref{eq:gen-case-mults}). Since the counting problem is completely solved when the ambient Lie algebra is finite-dimensional \cite[Tables 9-11]{dynkin}, this solves our original problem. In particular, the number of inequivalent \pisystems of type $K$ in $\lieg$ is necessarily  finite. 

Instances of our theorems encompass many examples of interest in physics, such as \pisystems of types $\ovex{A}_1 (=AE_3$), $\ovex{D}_8 (=DE_{10}) $ \cite{Henneaux2008} and $ \ovex{A}_n (=AE_n)$ \cite{klein-nico-aen} in $E_{10}$. We refer below to \S\ref{sec:physics-motivations} for a brief discussion of this connection.

\subsection{The physics of regular embeddings}\label{sec:physics-motivations}
\input{physics-embeddings}

\subsection{Brief historical remarks}
In Dynkin and Morita's original definitions, a \pisystem was required to be linearly independent (but see \cite[Table 7]{dynkin}, \cite[Chap 4, \S 2, exercises 29-38]{oni-vin}).  Oshima \cite{oshima} and Dynkin-Minchenko \cite{dyn-min} obtained extensions and variations of the results of \cite{dynkin}.

In the symmetrizable Kac-Moody context, Morita \cite{morita} and Naito \cite{naito} obtained the key initial results. A decade later, Feingold-Nicolai \cite{feingold-nicolai} rediscovered the definition of \pisystems, but imposed the restriction that all roots of a \pisystem be positive. They did not require linear independence, but as was pointed out by Henneaux et al \cite[\S 4.3]{HLPS}, their main theorem on embeddings arising out of \pisystems is false unless this condition is imposed. Our Theorem~\ref{thm:pisys-homo} is the corrected statement, in the more general setting of \pisystems that are not necessarily subsets of the positive real roots. Our Theorem~\ref{thm:pos-neg} serves as a link between the definitions of Morita and Feingold-Nicolai.

\subsection{Acknowledgements}
The authors are grateful to Thibault Damour for suggesting that we try to prove  conjugacy of embeddings of $\ovex{A}_1$ into $E_{10}$,  and for many illuminating discussions. We also wish to thank Ling Bao, Paul Cook and Axel Kleinschmidt for helpful discussions.

\section{\pisystems}\label{sec:pisys}

\subsection{} An integer matrix $A=(a_{ij})$ of size $n\times n$, where $n$ is a positive integer, is called a {\em generalized Cartan matrix\/}, {\em GCM\/} for short, if the following conditions are satisfied:
\begin{enumerate}
\item $a_{ii}=2$ for all $1\leq i\leq n$
\item $a_{ij}\leq 0$ whenever $1\leq i, j\leq n, i \neq j$
\item $a_{ij}=0$ if 
	$a_{ji}=0$ for $1\leq i, j\leq n$
\end{enumerate}
Given a GCM $A$ of size $n$,  we let $\algA$ denote the {\em Kac-Moody Lie algebra\/} associated to $A$ \cite[\S1.3]{kac}, with Cartan subalgebra $\lieh(A)$ and Chevalley generators $e_i, f_i$ for $1 \leq i \leq n$. We use terminology and notation as in the early chapters of~\cite{kac} without any further comment.

Let $\algprimeA$ denote the derived subalgebra $[\algA,\algA]$ of $\algA$. Let $\alpha_i(A), 1 \leq i \leq n$ denote the simple roots of $\algA$ and let $Q(A)$ be its root lattice, i.e., the free abelian group generated by the $\alpha_i(A)$. Both $\algA$ and $\algprimeA$ are $Q(A)$-graded Lie algebras, with $\deg e_i = \alpha_i(A) = - \deg f_i$ and $\deg h=0$ for all $h \in \lieh(A)$ \cite[Chapter 1]{kac}. We let $\Delta, \reroots, \imroots$ denote the sets of roots, real roots and imaginary roots respectively. For a root $\alpha$, we let $\algA_\alpha$ denote the corresponding root space. Each real root $\alpha$ defines a reflection $s_\alpha$ of $\lieh^*$ by $s_\alpha(\lambda) = \lambda - \aform{\alpha^\vee}{\lambda} \, \alpha$ where $\alpha^\vee \in \lieh(A)$ is the coroot corresponding to $\alpha$. The Weyl group $W(A)$ is the subgroup of $GL(\lieh^*)$ generated by the $s_\alpha, \; \alpha \in \reroots$.

\subsection{Multisets of real roots}\label{sec:msigma}
Let $A$ be a GCM, and let $\Sigma = \{\beta_1, \beta_2, \cdots, \beta_m\}$ be a finite sequence of real roots of $\algA$ (possibly with repetitions).
We define the $m \times m$ matrix
\[ M(\Sigma) := \left[ \aform{\beta_i^\vee}{\beta_j} \right]_{i,j=1}^m \]
We note that this is not a GCM in general. We let $\Sigma^\vee := \{\beta^\vee_1, \beta^\vee_2, \cdots, \beta^\vee_m\}$ be the corresponding multiset of coroots. Viewing these as real roots of $\lieg(A^T)$, we observe $M(\Sigma^\vee) = M(\Sigma)^T$.

A reordering of the elements of $\Sigma$ corresponds to a simultaneous permutation of the rows and columns of the matrix: $M(\Sigma) \mapsto P\,M(\Sigma)\,P^T$ for some $m \times m$ permutation matrix $P$. We will most often identify two such matrices without explicit mention.
\subsection{\pisystems}

\bdefn\mylabel{d:pisys}  Let $A$ be a GCM.
A {\em \pisystem} in $A$ is a finite collection of distinct real roots $\{\beta_i\}_{i=1}^m$ of $\lieg(A)$ such that $\beta_i-\beta_j$ is not a root for any $1\leq i\neq j\leq m$.
\edefn

This definition is essentially due to Dynkin \cite{dynkin} (for $A$ of finite type) and Morita \cite{morita} (in general), both of whom require that the $\{\beta_i\}_{i=1}^m$ be linearly independent; Morita calls such sets {\em fundamental subsets of roots}. The following proposition is stated in Morita (for the linearly independent case) without proof (see also Naito \cite{naito}). We supply the easy details.

\bprop\mylabel{p:bgcm}
Let $A$ be a GCM, and $\Sigma = \{\beta_i\}_{i=1}^m$ be a \pisystem in $A$.
Then the matrix $M(\Sigma)$ is a GCM.

\eprop
\bmyproof
	For any real root $\beta$ we have $\langle \beta^\vee,\beta\rangle=2$.  Indeed,  letting $\beta=w\alpha$ for a simple root~$\alpha$ and $w$~an element of the Weyl group,   we have $\beta^\vee=w(\alpha^\vee)$, and $\langle \beta^\vee,\beta\rangle=\langle w(\alpha^\vee),w\alpha\rangle=\langle \alpha^\vee,\alpha\rangle=2$.
	Suppose $\beta$ and $\gamma$ are distinct real roots such that $\gamma-\beta$ is not a root.  Consider $\{\gamma-p\beta, \ldots, \gamma+q\beta\}$ the ``$\beta$-string through~$\gamma$'' (\cite[Prop.~5.1]{kac}). Clearly $p=0$ and $\langle \beta^\vee,\gamma\rangle=p-q\leq0$.

	With $\beta$ and $\gamma$ as in the previous paragraph, if $\langle \beta^\vee,\gamma\rangle =0$, then $q=0$,   so that $\beta+\gamma$ is not a root, so the
	$\gamma$-string $\{\beta-p'\gamma,\ldots,\beta+q'\gamma\}$ through $\beta$ consists only of $\beta$, and so $\langle\gamma^\vee,\beta\rangle=p'-q'=0$.
\emyproof
\bdefn
We call $B:=M(\Sigma)$ the {\em type} of $\Sigma$, and refer to $\Sigma$ as a {\em \pisystem of type $B$ in $A$\/}.
\edefn

\subsection{Symmetrizable GCMs and \pisystems}\label{sec:symm-gcm}
An $n\times n$ GCM $A$ is {\em symmetrizable\/} if there exists a diagonal $n\times n$ matrix $D$ with positive rational diagonal entries such that $DA$ is symmetric. 
Let $\Sigma =\{\beta_i: 1 \leq i \leq m\}$ be a \pisystem of type $B$ in $A$. We note that if $A$ is a symmetrizable GCM, then so is $B$. Fix a choice of diagonal matrix $D$ which symmetrizes $A$, and let  $\form{\cdot}{ \cdot}$ denote the corresponding symmetric bilinear form  on $Q(A) \otimes_\integers \complex$, defined by:
\begin{equation} \label{eq:form-def} \form{\alpha_i(A)}{\alpha_j(A)} = D_{ii} \,a_{ij} \end{equation}
Since the $\beta_i$ are real roots of $\algA$, we know by \cite[Chapter 5]{kac} that:
\[b_{ij} = \langle \beta_i^\vee,\beta_j\rangle = \frac{2\form{\beta_i}{\beta_j}}{ \form{\beta_i}{\beta_i} }\]
Thus, $D' = \mathrm{diag}( \form{\beta_i}{\beta_i}/2)$ is a diagonal matrix with positive rational entries that symmetrizes $B$.  This choice of symmetrization defines a symmetric bilinear form on $Q(B)\otimes_\integers \complex$. As in equation~\eqref{eq:form-def} above, this is given by $\form{\alpha_i(B)}{\alpha_j(B)} = D'_{ii}\, b_{ij} = \form{\beta_i}{\beta_j}$. In other words, given the compatible choices of symmetrizations $(D,D')$ as above, the $\complex$-linear map
\begin{equation} \label{eq:def-q}
  Q(B) \otimes_\integers \complex \to Q(A) \otimes_\integers \complex, \;\;\;\;\alpha_i(B) \mapsto \beta_i \text{ for } 1 \leq i \leq m
\end{equation}
is form preserving.

\bdefn
We denote by $\qsig$ the map in equation~\eqref{eq:def-q} 
\edefn

Given $\alpha \in Q(A) \otimes_\integers \complex$ with $\form{\alpha}{\alpha} \neq 0$, the corresponding reflection $s_\alpha$ is given  by:
\[
s_\alpha(\gamma) = \gamma - \frac{2\form{\gamma}{\alpha}}{\form{\alpha}{\alpha}} \,\alpha \;\;\;\;\;\text{ for } \gamma \in Q(A) \otimes_\integers \complex.
\]
We note that $\qsig(s_\alpha(\beta)) = s_{\alpha'}(\beta')$ where $\alpha,\beta \in Q(B) \otimes_\integers \complex$ and $\alpha', \beta'$ are their images under $\qsig$.

\bthm\label{thm:pisys-homo}
Let $A$ be an $n\times n$ symmetrizable GCM and $\Sigma = \{\beta_i\}_{i=1}^m$ a \pisys\ of type $B$ in $A$. For each $1 \leq i \leq m$, let $e_{\beta_i}$, $e_{-\beta_i}$ be non-zero elements  in the root spaces $\algA_{\beta_i}$  and $\algA_{-\beta_i}$ respectively, such that $[e_{\beta_i},e_{-\beta_i}]=\beta_i^\vee$. Let $e_i, f_i$ be the Chevalley generators of $\algB$ and $h_i = [e_i,f_i]$. Then there exists a unique Lie algebra homomorphism $\isig:\algprimeB \rightarrow\algprimeA$ such that $e_i\mapsto e_{\beta_i}$, $f_i\mapsto e_{-\beta_i}$, $h_i \mapsto \beta_i^\vee$.
\ethm

\bmyproof
Since $A$ is symmetrizable, so is $B$, and $\algprimeB$ is generated by $e_i$, $f_i$, $h_i$, $1\leq i\leq m$ subject to the relations \cite[Theorem~9.11]{kac}:
\begin{gather}\label{e:defgb:1}
	[h_i,e_j]=b_{ij} \; e_j \quad\quad [h_i,f_j]=-b_{ij} \; f_j  
	\\
	\label{e:defgb:2}
	[h_i,h_j]=0 \\
	\label{e:defgb:3}
	[e_i,f_j]=\delta_{ij} \; h_i \quad\quad\quad \textup{and}\\
	\label{e:defgb:4}
	(\ad{e_i})^{1-b_{ij}} e_j= (\ad{f_i})^{1-b_{ij}} f_j = 0
\end{gather}
Any Lie algebra homomorphism from $\algprimeB$ is thus determined by the images of $e_i, f_i$ and $h_i$ ($1\leq i\leq m$).     Thus there is at most one Lie algebra homomorphism with the requisite properties.

To show that there exists such a homomorphism,  we need only verify that the relations in (\ref{e:defgb:1}) through (\ref{e:defgb:4}) are satisfied.
Relations (\ref{e:defgb:1}) and (\ref{e:defgb:2}) are clearly satisfied.    As for (\ref{e:defgb:3}) we consider two cases:   if $j=i$,  then it follows since $[e_{\beta_i},e_{-\beta_i}]=\beta_i^\vee$;   if $j\neq i$,  then it follows since $\beta_i-\beta_j$ is not a root of $\algA$ by the definition of \pisystem.   As for (\ref{e:defgb:4}),   it follows from the fact \cite[Prop.~5.1]{kac} that the $\beta_i$-string through $\beta_j$ consists of $\beta_j$, $\beta_j+\beta_i$, \ldots, $\beta_j+k\beta_i$,  where $k=\langle\beta_i^\vee, \beta_j\rangle$.
\emyproof

The following proposition is equivalent to that of Naito \cite[Theorem 3.6]{naito}, though his proof is different (without using the Serre relations). In the interest of completeness, we give a (slightly simpler) argument.
\begin{proposition} \label{prop:pisys-li}
With notation as in the above theorem, if $\Sigma$ is linearly independent (in $Q(A) \otimes_\integers \complex$), one can extend the map $\isig$ to an injective map from $\algB$ to $\algA$. 
\end{proposition}

\begin{proof}
Suppose that $\{\lieh; \alpha_1^\vee,\ldots,\alpha_n^\vee;\alpha_1,\ldots,\alpha_n\}$ is a realization of~$A$ \cite[Chapter 1]{kac}.   Let $\liek$ be any subspace of $\lieh$ of smallest possible dimension such that
			(i) $\liek$ contains $\beta_1^\vee$, \ldots, $\beta_m^\vee$, and
			(ii) the restrictions of $\beta_1$, \ldots, $\beta_m$ to $\liek$ are linearly independent as elements of~$\liek^*$
 (this is possible since we are given that the $\beta_i$ are linearly independent). Then
  \be
  \item\label{i:one:p:rankb} $(\liek, \beta_1^\vee,\ldots,\beta_m^\vee;\beta_1|_\liek,\ldots,\beta_m|_\liek)$ is a realization of $B$.
				\item\label{i:two:p:rankb} $\rank{B}\geq \rank{A}-2(n-m)$.
  \ee

  Assertion~(\ref{i:one:p:rankb}) follows easily from the definition of realization.  As for assertion~(\ref{i:two:p:rankb}), observe that $\{\beta_i^\vee\}_{i=1}^m$ is in the span of $\{\alpha_i^\vee\}_{i=1}^n$:  this follows from the definition of $\beta^\vee$ for a real root $\beta$ as $w(\alpha_i^\vee)$ where $w$ is an element of the Weyl group such that $\beta=w(\alpha_i)$.  We have $B=YAX$, where $X=(x_{ij})$ is the $n\times m$ matrix such that $\beta_j=\sum_{i=1}^n x_{ij}\alpha_i$ and $Y=(y_{ij})$ is the $m\times n$ matrix such that $\beta_j^\vee=\sum_{i=1}^ny_{ji}\alpha_i^\vee$. The matrices $X$ and $Y$ are both of rank~$m$.   The assertion now follows easily from elementary linear algebra.

Now, $\lieg(B)$ is generated by $\liek$, $e_i, f_i$ subject to the relations specified in the proof of Theorem~\ref{thm:pisys-homo}
together with the following:
\[ [k,e_i]=\beta_i(k)e_i \quad\quad\quad [k,f_i]=-\beta_i(k)f_i  \quad\quad\quad\quad 
	[k_1,k_2]=0 \quad\quad\quad\textup{for $k, k_1, k_2$ in $\liek$}\]
We map $\liek$ to $\lieh$ via the natural inclusion; $e_i, f_i$ are mapped to $e_{\beta_i}, e_{-\beta_i}$ as before. We only need to check that the additional relations above hold. But these are obvious.

Finally,  we show that the homomorphism is an embedding.   The kernel of the homomorphism being an ideal of $\algB$,   it either contains the derived algebra $\deralgB$ or is contained in the center~\cite[\S1.7(b)]{kac}.     Since $e_i\mapsto e_{\beta_i}$ (and $e_i$ is contained in $\deralgB$ by (\ref{e:defgb:1})) the first possibility is ruled out.    Thus the kernel is contained in the center.    But the center is contained in the subspace $\liek$ (\cite[Prop.~1.6]{kac}) and on $\liek$ the homomorphism is an inclusion.  Thus the kernel is zero.
\end{proof}

\bremark \label{rem:easy} The following easy observations are often useful:
  \be
\item  If $\Sigma$ is linearly independent, then $\qsig$ is an injection.
  \item If $\det B \neq 0$, then $\Sigma$ is linearly independent.

    \ee
\eremark

\bexamplenobox \label{eg:pisys-examples}
  \be
\item[(i)] Let $A$ be a GCM of finite type. Dynkin \cite{dynkin} showed that if $\mathfrak{m}$ is a regular semisimple subalgebra of $\lieg(A)$, then there exists a GCM $B$ of finite type and a \pisystem $\Sigma$ of type $B$ in $A$ such that $\mathfrak{m} = \isig(\algB)$.

  \smallskip
\item[(ii)] Let us take $A=[2]$, so that $\algA = \algprimeA = \mathfrak{sl}_2 \complex$. Let $\Sigma = \{\alpha_1, -\alpha_1\} = \roots(A)$. This is clearly a \pisystem in $A$, of type $B = \begin{bmatrix} 2 & -2 \\ -2 & 2 \end{bmatrix}$. The corresponding Kac-Moody algebra $\algB$ is the affine Lie algebra $\widehat{\mathfrak{sl}_2 \complex}$. We then have \cite[Chapter 7]{kac}, $\algprimeB = \mathfrak{sl}_2 \complex \,\otimes\, \complex[t,t^{-1}] \oplus \complex c$, the universal central extension of the loop algebra of $\mathfrak{sl}_2$. The generators of $\algprimeB$ are $ e_1 = X, f_1 = Y, e_2 = Y \otimes t, f_2 = X \otimes t^{-1}$, where $X=\left(\begin{smallmatrix} 0 & 1 \\ 0 & 0 \end{smallmatrix}\right)$ and $Y=\left(\begin{smallmatrix} 0 & 0 \\ 1 & 0 \end{smallmatrix}\right)$ are the standard generators of $\mathfrak{sl}_2 \complex$.
    The map defined in Theorem~\ref{thm:pisys-homo} may be chosen to be: 
    \[e_1 \mapsto X, \, f_1 \mapsto Y, \, e_2 \mapsto Y, \, f_2 \mapsto X\]

  \item[(iii)] More generally, let $A$ be any finite type GCM and $\algA$ the corresponding finite dimensional simple Lie algebra, with highest root $\theta$.
    Consider the \pisystem $\Sigma$ consisting of the simple roots of $\algA$ together with $-\theta$. This has type $B$, the GCM of the untwisted affinization of $\algA$. The map $\algprimeB = \algA \otimes \complex[t,t^{-1}] \, \oplus \complex c \to \algA$ defined by Theorem~\ref{thm:pisys-homo} may be chosen to be the {\em evaluation map} at $t=1$:
    \[c \mapsto 0 \text{ and } \zeta \otimes g(t) \mapsto g(1) \, \zeta \text{ for all } \zeta \in \algA, \,g \in \complex[t,t^{-1}]. \]
\ee
\eexamplenobox

\begin{lemma}\label{lem:ideal-rootvecs}
  Let A be an $n\times n$ GCM. Let $I$ be an ideal of  $\algprimeA$ that does not contain any of the Chevalley generators i.e., $e_i, f_i \not\in I$ for all $i$. Then $\algprimeA_\alpha \cap I = (0)$ for all roots $\alpha$.
\end{lemma}

\begin{proof}

Suppose $\alpha$ is a positive, non-simple root. Assume $e_{\alpha}\in I$ for some  nonzero $e_\alpha \in \algprimeA_\alpha$. By \cite[Lemma 1.5]{kac}, there exists $i_1$ such that $[f_{i_1},e_{\alpha}]\neq 0$. If $\alpha-\alpha_{i_1}$ is not a simple root, we can find $i_2$ such that $[f_{i_2},[f_{i_1},e_{\alpha}]]\neq 0$. Proceeding this way, after finitely many steps we get $[f_{i_k}[\cdots f_{i_2},[f_{i_1},e_{\alpha}]]\cdots]=e_i\in I$, which contradicts the hypothesis on $I$. If $\alpha$ were a negative root to begin with, the proof is analogous.
\end{proof}

\bremark
\begin{enumerate}
\item Let $I$ be an ideal of $\algprimeA$. We observe that if $I$ contains one of $e_i, f_i, \alpha_i^\vee$, then it contains all three. 
\item If $A$ is an indecomposable GCM, then any proper ideal of $\algprimeA$ satisfies the hypothesis of lemma~\ref{lem:ideal-rootvecs}. To see this, suppose $e_i$ is in $I$. Then, so are $f_i$ and  $\alpha_i^{\vee}$. Since $A$ is indecomposable, for each fixed $j$, there exist $i_1,i_2,\cdots i_s$ such that $a_{i i_1}a_{i_1 i_2}\cdots a_{i_s j}\neq 0$. Since, $[\alpha_i^{\vee},e_{i_1}]=a_{ii_1}e_{i_1}$, we conclude $I$ contains $e_{i_1}$, and hence also $f_{i_1},  \alpha_{i_1}^{\vee}$. Proceeding in this manner, we get $e_j,f_j,\alpha_j^{\vee} \in I$. Since this holds for all $j$, we obtain $I=\algprimeA$, a contradiction.
\end{enumerate}
\eremark

While the map $\isig$ of Theorem \ref{thm:pisys-homo} need not be injective when $\Sigma$ is linearly dependent, we nevertheless have the following useful result which states that it is injective on each root space.
\begin{corollary}\label{cor:inj_root_spaces}
  The map $\isig: \,\algprimeB \rightarrow\algprimeA$ defined in Theorem \ref{thm:pisys-homo} is injective when restricted to $\algprimeB_\alpha$ for $\alpha \in \Delta(B)$.
Further, the image of $\algprimeB_\alpha$ is contained in $\algprimeA_{\qsig(\alpha)}.$
\end{corollary}

\begin{corollary} \label{cor:real-to-real}
\be
\item $\qsig(\reroots(B)) \subset \reroots(A)$ and $\qsig(\imroots(B)) \subset \imroots(A) \cup \{0\}$.
\item If further $\Sigma$ is linearly independent, then  $\qsig(\imroots(B)) \subset \imroots(A)$.
\ee
\end{corollary}
\begin{proof}
  Corollary \ref{cor:inj_root_spaces} implies that if $\alpha$ is a root of $\algprimeB$, then $\qsig(\alpha)$ is either 0 or a root of $\algprimeA$. Since $\form{\alpha}{\alpha} = \form{\qsig(\alpha)}{\qsig(\alpha)}$ and real roots are precisely the roots of positive norm \cite[Proposition 5.2]{kac}, we conclude that real roots map to real roots and imaginary roots to imaginary roots or 0. The second part is now obvious from Remark~\ref{rem:easy}.
\end{proof}  

The above corollary, for linearly independent $\Sigma$ was first obtained by Naito \cite[Theorem 3.8]{naito}.
Next, we have the converse to Theorem~\ref{thm:pisys-homo}:

\begin{proposition}\label{prop:converse-pisys-homo}
Let $A_{n\times n}, B_{m \times m}$ be symmetrizable GCMs. Let $e_i, f_i$ denote the Chevalley generators of $\algB$. Suppose $\phi: \algprimeB \to \algprimeA$ is a Lie algebra homomorphism satisfying $0 \neq \phi(e_i) \in \algprimeA_{\beta_i}$,  $0 \neq \phi(f_i) \in \algprimeA_{-\beta_i}$ for all $1 \leq i \leq m$, for some real roots $\{\beta_i\}_{i=1}^m$ of $\algprimeA$. Then, the set $\Sigma=\{\beta_i\}_{i=1}^m$ is a \pisystem of type $B$ in $A$.
\end{proposition}
\begin{proof}
  Given a real root $\beta$ and any root $\gamma$ of $\algprimeA$, it follows from elementary $\mathfrak{sl}_2$ theory (applied to the $\beta$-string through $\gamma$) that
  \begin{equation} \label{eq:key-fact}
    [\algprimeA_\beta, \,\algprimeA_\gamma] \neq 0 \text{ iff } \beta + \gamma  \text{ is a root of } \algprimeA
   \end{equation}

  Now, since $[e_i,f_j]=0$ for $1 \leq i \neq j \leq m$, we apply $\phi$ to conclude that $[\algprimeA_{\beta_i}, \algprimeA_{-\beta_j}] =0$. Hence $\beta_i - \beta_j$ is not a root of $\algprimeA$, and $\Sigma$ is thus a \pisystem.

\smallskip
Next, we show that the type of this \pisystem is exactly $B$. Note that $|\langle\beta_i^\vee, \beta_j\rangle|$ is the largest integer $k$ for which $\beta_j + k' \beta_i$ is a root of $\algprimeA$ for $0 \leq k' \leq k$.
Let $\alpha_i(B)$ denote the simple roots of $\algprimeB$; their images under $\qsig$ are the $\beta_i$. We have $\ell = |b_{ij}|$ is the largest integer for which $\alpha_j(B) +  \,\ell'\alpha_i(B)$ is a root of $\algprimeB$ for $0 \leq \ell' \leq \ell$. In fact $\gamma = \alpha_j(B) +  \,\ell\alpha_i(B) \in \reroots(B)$, and by corollary~\ref{cor:real-to-real}, $\qsig(\gamma)\in \reroots(A)$. Thus, $k \geq \ell$.

\smallskip
By \eqref{eq:key-fact} above, $[\algprimeB_{\alpha_i(B)}, \,\algprimeB_\gamma] =0$, and since these two real root spaces map isomorphically to the corresponding real root spaces of $\algprimeA$, we conclude $[\algprimeA_{\beta_i}, \,\algprimeB_{\qsig(\gamma)}] =0$. By \eqref{eq:key-fact} again, $\beta_i + \qsig(\gamma) = \beta_j + \,(\ell+1) \beta_i $ is not a root of $\algprimeA$. Hence $k \leq \ell$, and we obtain $\langle\beta_i^\vee, \beta_j\rangle = b_{ij}$ as required.
\end{proof}

\begin{corollary}\label{cor:transitivity}
Let $A, B, C$ be symmetrizable GCMs. If $A$ has a \pisys of type $B$ and $B$ has a \pisys of type $C$,  then $A$ has a \pisys of type $C$.
\end{corollary}
\bmyproof
Theorem~\ref{thm:pisys-homo} gives us Lie algebra morphisms $\algprimeA[C] \to \algprimeB \to \algprimeA$. By corollary~\ref{cor:inj_root_spaces}, both these maps are injective on real root spaces. The generators $e_i, f_i$ of $\algprimeA[C]$ map to real root vectors of $\algprimeB$. Thus, under the composition of these two morphisms, $e_i, f_i$  map to non-zero real root vectors of $\algprimeA$. The corresponding roots are clearly negatives of each other. Proposition~\ref{prop:converse-pisys-homo} now completes the proof. 
\emyproof

\smallskip
If $\Sigma_1, \Sigma_2$ denote the \pisystems of the above corollary, of types $B$ and $C$ respectively, then the \pisystem of type $C$ in $A$ that one obtains from the proof above is just $\qsig[\Sigma_1](\Sigma_2)$.

\subsection{}\label{sec:regsub-rmk}
As mentioned in the introduction, \pisystems were first defined by Dynkin in his study of regular semisimple subalgebras of semisimple Lie algebras. In this setting, any set of simple roots of a closed subroot system of the root system (of a semisimple Lie algebra) is a \pisystem. The converse is also true, as can be seen from Theorem \ref{thm:pisys-homo}.

In the infinite dimensional setting, Naito \cite{naito} defined a {\em regular subalgebra} of a Kac-Moody algebra $\lieg(A)$ to be any subalgebra of the form $\widehat{i}_\Sigma(\lieg(B))$ where $\Sigma$ is a linearly independent \pisystem of type $B$ in $A$ and $\widehat{i}_\Sigma$ denotes an extension of $i_\Sigma$ as in Proposition~\ref{prop:pisys-li}.

\section{Weyl group action on \pisystems}
\subsection{}
We freely use the notation of \S\ref{sec:pisys}. Let $A$ be a symmetrizable GCM. Let $W(A)$ denote the Weyl group of $A$. It acts on the set of roots of $A$, preserving each of the subsets of real and imaginary roots. Further this action preserves the symmetric bilinear form $\form{\cdot}{\cdot}$ on $\lieh^*(A)$. Thus, there is an induced action of $W(A)$  on the set of all \pisystems\ in~$A$ of a given type $B$.

\subsection{}
When $A$ is of finite type, it is easy to see that every linearly independent \pisystem in $A$ is $W(A)$-conjugate to a \pisystem contained in
the set of {\em positive roots} of $A$. To see this, take an element $\gamma \in \lieh^*(A)$ such that $\form{\gamma}{\alpha} >0$ for all elements $\alpha$ of the \pisystem. The element $w \in W(A)$ which maps $\gamma$ into the dominant Weyl chamber will clearly also map the \pisystem to a subset of the positive roots.

This proof fails in the general case; such $w$ does not exist unless $\gamma$ is in the Tits cone. For instance, the negative simple roots of $A$ form a \pisystem of type $A$ in $A$. This set cannot be $W(A)$-conjugated to a subset of positive roots if $A$ is not of finite type; this can be seen using for instance \cite[Theorem 3.12c]{kac}.
The next theorem shows that this is essentially the only obstruction. 
\bthm \label{thm:pos-neg}
Let $A, B$ be symmetrizable GCMs and $\Sigma$ a linearly independent $\pi$-system of type $B$ in $A$. If $B$ is {\em indecomposable}, then:
\be
\item There exists $w\in W(A)$ such that $w\Sigma \subset  \reroots[+](A)$ or $w\Sigma \subset \reroots[-](A)$.
\item There exist $w_1, w_2 \in W(A)$ such that $w_1 \Sigma \subset \reroots[+](A)$ and $w_2\Sigma \subset \reroots[-](A)$ if and only if $B$ is  of finite type.
\ee
\ethm

The proof occupies the next two subsections.
\subsection{}
The proof of theorem~\ref{thm:pos-neg} closely follows that of \cite[Proposition 5.9]{kac}. The first part of this theorem, in the special case $|\Sigma|=2$ was proved by Naito in \cite{naito}. We first recall some relevant facts about the roots of a Kac-Moody algebra. Let $B$ be an {\em indecomposable} GCM, and let $\algB$ denote the corresponding Kac-Moody algebra.  Let $Q(B)$ denote its root lattice. We use the notation introduced already for the sets of roots, real roots, positive roots etc. Let $\real_+$ denote the set of non-negative reals. Define:
\[\cim = \bigcup_{\alpha \in \imroots[+](B)} \real_+\alpha, \;\;\;\; \cre = \bigcup_{\alpha \in \reroots[+](B)} \real_+\alpha.\] 
We then have the following result due to Kac \cite[Proposition 1.8]{kac1978},  \cite[\S 5.8]{kac}:
\bprop\label{prop:kac-cim} (Kac) In the metric topology on the real span of $Q(B), \;\;\cimcl$ is the convex hull of the set of limit points of $\cre$.
In particular, it is a convex cone.
\eprop

Now suppose $Q(B) \subset \euc$ for some real vector space $\euc$. Let $\{\epsilon_i\}_{i=1}^n$ be a basis of $\euc$. Define $\euc[+]$ to be the $\real_+$ span of the $\epsilon_i$, and let $\euc[-] = -\euc[+]$.

\begin{lemma} \label{lem:imroots-pos-neg}
If $\roots(B) \subset \euc[+] \cup \euc[-]$, then $\imroots[+](B) \subset \euc[+]$ or  $\imroots[+](B) \subset \euc[-]$.
\end{lemma}
\begin{proof}
  Consider the set $\cimcl$; it has the following properties: (i) It is convex, by Proposition~\ref{prop:kac-cim}. (ii) It is contained in $\euc[+] \cup \euc[-]$, by the given hypothesis. (iii) It does not contain a line  (i.e., for nonzero $x \in \euc$, both $x$ and $-x$ cannot belong to this set), because $\cimcl \subset \real_+\left( \roots[+](B)\right)$.

  It is easy to see that these properties imply that $\cimcl$ must be entirely contained either in $\euc[+]$ or in $\euc[-]$.
  \end{proof}

Under the same hypothesis as lemma~\ref{lem:imroots-pos-neg}, we have:
\begin{lemma} \label{lem:reroots-fin-nonpos}
If $\imroots[+](B) \subset \euc[+]$, then all but finitely many real roots of $B$ lie in $\euc[+]$.
\end{lemma}
\begin{proof}
  First, we define an inner product on $\euc$ by requiring the $\epsilon_i$ to be an orthonormal basis. This defines the standard metric topology on $\euc$, and thereby on the $\real$-span of $Q(B)$.
  
  Let $M := \reroots[+](B) \cap \euc[-]$, and $\widehat{M} := \{\alpha/\norm{\alpha}: \alpha \in M\}$. Here, the norm is that of the Euclidean space $\euc$. Observe that $\widehat{M}$ is a subset of $\cre \cap \euc[-] \cap S$, where $S$ is the unit sphere in $\euc$. If $\widehat{M}$ is an infinite set, then, it has a limit point, say $\zeta$. Now $\zeta \in \euc[-] \cap S$, and by Proposition~\ref{prop:kac-cim}, $\zeta \in \cimcl$. But $\cimcl \subset \euc[+]$ by hypothesis. This contradiction establishes the lemma. 
\end{proof}
\begin{proposition} \label{prop:w-conj-pos-neg}
Let $\roots(B) \subset \euc[+] \cup \euc[-]$. There exists $w \in W(B)$ such that $w\roots[+](B) \subset \euc[+]$ or $w\roots[+](B) \subset \euc[-]$.
\end{proposition}
\begin{proof}
  By lemma~\ref{lem:imroots-pos-neg}, the positive imaginary roots are all contained in $\euc[+]$ or in $\euc[-]$; we may suppose (replacing the $\epsilon_i$ with their negatives if need be) that   $\imroots[+](B) \subset \euc[+]$. Consider $F := \reroots[+](B) \cap \euc[-]$; this is finite  by lemma~\ref{lem:reroots-fin-nonpos}. If this set is non-empty, it contains some simple root $\alpha$ of $\algB$. Since the simple reflection $s_\alpha$ defines a bijective self-map of $\reroots[+](B) \backslash \{\alpha\}$, it is clear that $F' := s_\alpha\left(\reroots[+](B)\right) \cap \euc[-]$ contains one fewer element than $F$. Iterating this procedure, we can find $w$, a product of simple reflections, such that $ w\,\reroots[+](B) \cap \euc[-]$ is empty, as required.
\end{proof}

\subsection{}
Finally, we are in a position to prove theorem~\ref{thm:pos-neg}. With notation as in the theorem, observe that the linear independence of $\Sigma$ implies that $\qsig: Q(B) \to Q(A)$ is injective. By corollary~\ref{cor:real-to-real}, $\qsig(\roots(B)) \subset \roots(A) = \roots[+](A) \cup \roots[-](A)$. We define $\euc$ to be the $\real$-span of $\roots(A)$ and take $\{\epsilon_i\}$ to be the basis of simple roots of $\algA$. Then, clearly, $\qsig(\roots(B)) \subset \euc[+] \cup \euc[-]$. Identifying $\roots(B)$ with its image under $\qsig$, and appealing to proposition~\ref{prop:w-conj-pos-neg} completes the proof of part (1).

To prove part (2), since $w_1 \Sigma \subset \reroots[+](A)$, we have $w_1 (\qsig (\roots[+](B))) \subset \roots[+](A)$. Consider the set $R:=\qsig(\imroots[+](B))$. We have (i) $R \subset \imroots(A)$, by corollary~\ref{cor:real-to-real}, and (ii) $w_1R \subset \roots[+](A)$. Since the sets $\imroots[\pm](A)$ are both $W(A)$-invariant, this implies $R \subset \imroots[+](A)$.
Similarly, from $w_2\Sigma \subset \reroots[-](A)$, we conclude $R \subset \imroots[-](A)$. This means $R$ is empty, or in other words, that $B$ is of finite type.

Conversely, if $B$ is of finite type, then $\roots[+](B)$ is finite. Hence its intersections with  $\roots[+](A)$ and $\roots[-](A)$ are both finite sets. The proof of Proposition~\ref{prop:w-conj-pos-neg} shows that there exist elements of $W(A)$ which map $\roots[+](B)$ to subsets of $\roots[\pm](A)$.
\qed

\subsection{}
As is evident from Example~\ref{eg:pisys-examples}(ii), the conclusion of Theorem~\ref{thm:pos-neg} is false if $\Sigma$ is not assumed to be linearly independent, even when $A$ is of finite type.

\subsection{}  Let $A, B$ be symmetrizable GCMs. A \pisystem $\Sigma$ of type $B$ in $A$ is said to be {\em positive} (resp. {\em negative}) if it is $W(A)$-conjugate to a \pisystem all of whose elements are positive (respectively negative) roots. Theorem \ref{thm:pos-neg} implies that if $\Sigma$ is linearly independent and $B$ is indecomposable and not of finite type, then $\Sigma$ is either positive or negative, but not both. We record below a simple criterion to determine the sign that was obtained in the course of the proof of Theorem~\ref{thm:pos-neg}.
\begin{proposition} \label{prop:criterion-pos-neg}
Let $A, B$ be symmetrizable GCMs, with $B$ indecomposable and not of finite type. Let $\Sigma$ be a linearly independent \pisystem of type $B$ in $A$. Then the following are equivalent:
\be
\item $\Sigma$ is positive (resp. negative).
\item $\qsig(\alpha) \in \imroots[+](A)$ (resp. $\imroots[-](A)$) for {\em every} $\alpha \in \imroots[+](B)$.
\item $\qsig(\alpha) \in \imroots[+](A)$ (resp. $\imroots[-](A)$) for some $\alpha \in \imroots[+](B)$.
\ee \qed
\end{proposition}

\subsection{}
Let $\mult(B,A)$ denote the number of $W(A)$-orbits of \pisystems of type $B$  in $A$ (this could be infinity in general).
When $A, B$ are of finite type, Borel-de Siebenthal and Dynkin determined the pairs for which $\mult(B,A)>0$. Dynkin went further, and also determined the values of $\mult(B,A)$; these turn out to be 1 for almost all cases, except for a few where it is 2 \cite[Tables 9-11]{dynkin}

\section{$\pi$-systems of affine type}
\subsection{} \label{sec:subdig-po}
Let $S(A)$ denote the Dynkin diagram associated to the GCM $A$ \cite[\S 4.7]{kac}. Any subset of the vertices of $S(A)$ together with the edges between them will be called a {\em subdiagram} of $S(A)$ (and we will use $\subseteq$ to denote the relation of being a subdiagram). Given $\alpha = \sum_{i=1}^n c_i \alpha_i$, we define $\supp \alpha$ to be the set $\{i: c_i \neq 0\}$ and view it as a subset of the vertices of $S(A)$. Given a subdiagram $Y$ of $S(A)$, we say $\alpha$ is supported in $Y$ if $\supp \alpha$ is contained in the set of vertices of $Y$. We also let $Y^\perp$ denote the subset of vertices of $S(A)$ that are not in $Y$ and are not connected by an edge to any vertex of $Y$.
If $\alpha$ is supported in $Y$ and $\beta$ in $Y^\perp$, then clearly $\form{\alpha}{\beta}=0$.
\begin{lemma}\label{lem:nullroot-perp}
Let $A$ be a symmetrizable GCM and $Y$ a subdiagram of $S(A)$ of affine type. Let $\delta_Y$ denote the null root of $Y$. If $\beta \in \roots(A)$ is such that $\form{\beta}{\delta_Y}=0$, then $\supp \beta \subset Y \sqcup Y^\perp$.
\end{lemma} 
\begin{proof}
We write $\beta = \sum_{p \in S(A)} c_p \alpha_p$, where the coefficients are all non-negative, or all non-positive. Let $\supp \beta$ denote the set of $p$ for which $c_p$ is nonzero. Now, $\form{\alpha_p}{\delta_Y}$ is 0 for $p \in Y$, and $\leq 0$ when $p \not\in Y$. Since all coefficients are of the same sign, every $p \in \supp \beta$ must be either in $Y$ or in $Y^\perp$.
\end{proof}

\bthm \label{thm:affine-pi}
Let $A$ be a symmetrizable GCM and $B$ be a GCM of affine type. Suppose $\Sigma$ is a linearly independent $\pi$-system of type $B$ in $A$. Then,
\be
\item There exists an affine subdiagram $Y$ of $S(A)$ and $w \in W(A)$ such that every element of $w\Sigma$ is supported in $Y$.
\item Suppose $(Y^\prime, w^\prime)$ is another such pair, i.e., with $Y^\prime$ a subdiagram of affine type, $w^\prime \in W(A)$ such that $w^\prime\Sigma$ is supported in $Y^\prime$. Then $Y = Y^\prime$ and $w^\prime w^{-1} \in W(Y \sqcup Y^\perp)$. 
\item $\mult(B,A) = \infty$.
  \ee
\ethm
\begin{proof}
 Let $\Sigma=\{\beta_i\}_{1}^m$. Let $\{\alpha_i(B)\}_{1}^m$ denote the simple roots of $\algB$ and let $\delta_B$ denote its null root. Let 
 $\delta_{\Sigma} = \qsig(\delta_B)$. By corollary~\ref{cor:real-to-real}(2) and the fact that $\qsig$ preserves forms, we obtain that $\delta_\Sigma$ is an isotropic root of $\algA$. By \cite[Proposition 5.7]{kac}, there exists $w\in W(A)$ such that $w(\delta_{\Sigma})$ is supported on an affine subdiagram $Y$ of $S(A)$ and $w(\delta_{\Sigma})=k\delta_Y$ for some nonzero integer $k$, where $\delta_Y$ is the null root of $Y$.

 Now, $0 = \form{\alpha_i(B)}{\delta_B} = \form{\beta_i}{\delta_\Sigma} = k\form{w\beta_i}{\delta_Y}$ for all $i=1, \cdots, m$. We conclude $\supp w\beta_i \subset Y\sqcup Y^\perp$, by lemma~\ref{lem:nullroot-perp}. Since $w\beta_i$ is a root, its support is connected, and hence contained entirely in $Y$ or entirely in $Y^\perp$. However, $w\Sigma$ is a \pisystem of type $B$, an indecomposable GCM. So, $w\Sigma$ cannot be written as a disjoint union of two mutually orthogonal subsets. This means that either $\supp w\beta_i \subset Y$ for all $i$, or  $\supp w\beta_i \subset Y^\perp$ for all $i$. The latter is impossible since $k\delta_Y= w(\delta_{\Sigma})$ is a positive integral combination of the $w\beta_i$.
 This proves part (1).

 Now, if $(Y',w')$ is another such pair, then since the only isotropic roots of $\algA$ supported on subsets of $Y'$ are the multiples of $\delta_{Y'}$, we obtain $w'(\delta_\Sigma) = k'\delta_{Y'}$ for $k' \neq 0$. Define $\sigma = w'w^{-1}$, so $\sigma(k\delta_Y) = k'\,\delta_{Y'}$.
Since $\delta_Y$ is a positive imaginary root of $\algA$, so is $\sigma\delta_Y$; thus $k$ and $k'$ have the same sign. We may suppose $k, k' >0$. Now $k\delta_Y$ and $k'\delta_{Y'}$ are antidominant weights (i.e., their negatives are dominant weights) of $\algA$, which are $W(A)$-conjugate. By \cite[Proposition 5.2b]{kac}, we get $k\delta_{Y}=k'\delta_{Y'}$. Thus,  $Y=Y'$, $k=k'$ and $\sigma\delta_Y = \delta_Y$.

Since $\delta_Y$ is antidominant, the simple reflections that fix $\delta_Y$ generate the stabilizer of $\delta_Y$ \cite[Proposition 3.12a]{kac}. By lemma~\ref{lem:nullroot-perp}, this stabilizer is just $W(Y \sqcup Y^\perp)$. Thus $\sigma \in W(Y \sqcup Y^\perp)$, proving part (2).

  Finally, let $\Sigma = w\Sigma$ denote the \pisystem  of part (1). Now $Y$ is of affine type, untwisted or twisted. In either case, from the description of the real roots of an affine Kac-Moody algebra \cite[Chap 6]{kac}, the following holds: $\reroots(Y) + 6p\, \delta_Y \subset \reroots(Y)$ for all $p \in \integers$. Consider
  \[\Sigma_p := \{\alpha + 6p\,\delta_Y: \alpha \in \Sigma\} \text{ for } p \in \integers. \]
  Since $\delta_Y$ is orthogonal to every root of $\lieg(Y)$, it is clear that $\Sigma_p$ is a linearly independent \pisystem of type $B$ in $A$, supported in $Y$.
  From the proof of part (1), we know $\qsig(\delta_B) = k \delta_Y$ for some nonzero integer $k$. From the definition of $\Sigma_p$, we obtain  \begin{equation} \label{eq:delta-reln}
    \qsig[\Sigma_p](\delta_B) = (k + 6ph) \delta_Y
    \end{equation}
where  $h$ is the Coxeter number of the affine Kac-Moody algebra $\algB$.
  We claim that the $\Sigma_p$ are pairwise $W(A)$-inequivalent. Suppose  $\Sigma_m$ and $\Sigma_n$ are in the same $W(A)$-orbit. Then, from part (2), we obtain $\Sigma_m = \sigma ( \Sigma_n) $ for some $\sigma \in W(Y \sqcup Y^\perp)$. In particular, this means $\qsig[\Sigma_m](\delta_B) = \sigma(\qsig[\Sigma_n](\delta_B))$. Since $\sigma$ fixes $\delta_Y$, equation~\eqref{eq:delta-reln} implies $m=n$. This completes the proof of part (3).
\end{proof}

\begin{corollary}
Let $A$ be a symmetrizable GCM such that $S(A)$ has no subdiagrams of affine type. Then $A$ contains no linearly independent \pisystems of affine type. 
\end{corollary}

This follows immediately from the proposition. We remark that  Figure~\ref{fig:irreg-hyps} contains examples of such $S(A)$.

\bremark
\be
\item The conclusion of theorem~\ref{thm:affine-pi} is false without the linear independence assumption, as in Example~\ref{eg:pisys-examples} (ii), (iii).
\item Let $A, B$ be symmetrizable GCMs, with $B$ of affine type. Suppose $A$ contains a linearly independent \pisystem of type $B$. Theorem~\ref{thm:affine-pi} implies that some affine type subdiagram $Y$ of $S(A)$ also contains a linearly independent \pisystem of type $B$. This allows us to determine the possible set of such $B$ in two steps: (i) find all affine subdiagrams $Y$ of $S(A)$, and (ii) for each such $Y$, list out all the $B$'s which occur as GCMs of linearly independent \pisystems of $Y$.
\item We note that step (ii) above can in-principle be carried out using the results of \cite{roy-venkatesh} (see also \cite{naito,FRT,dyer-lehrer}).
\ee
\eremark

\section{Hyperbolics and Overextensions}

\subsection{} Let $A$ be a symmetrizable GCM and $X=S(A)$ be its Dynkin diagram \cite[\S 4.7]{kac}. If $A$ is symmetric, we will call $X$ {\em simply-laced}.
\begin{definition}\label{ovex-def}
Let $Z$ be a simply-laced Dynkin diagram. We say that $Z$ is an {\em overextension} or of $\Ext$ type if there exists a vertex $p$ in $Z$ such that the subdiagram $Y = Z \backslash \{p\}$ is of affine type and $\form{\delta_Y}{\alpha_p} = -1$.
\end{definition}

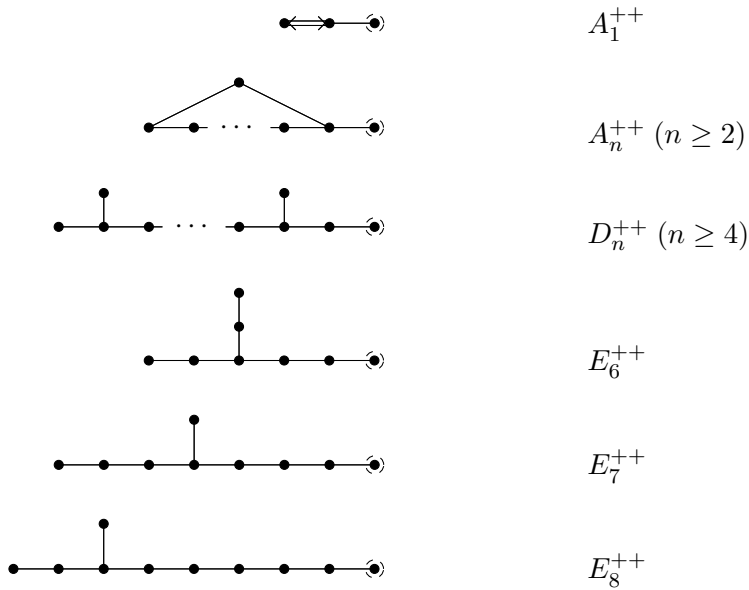
\begin{figure}[h]
\begin{tabular}{rp{2cm}l}
  \begin{tikzpicture}[scale=0.6]
    \filldraw[] (0,0) circle (.1);
    \filldraw[] (1,0) circle (.1);
     \filldraw[] (2,0) circle (.1);
      \filldraw[] (0,-.05) -- +(1,0);
      \filldraw[] (0,.05) -- +(1,0);
      \filldraw[] (1,0) -- +(1,0);
      \filldraw[] (.1,0) -- +(.15,.15);
      \filldraw[] (.1,0) -- +(.15,-.15);
      \filldraw[] (.9,0) -- +(-.15,.15);
      \filldraw[] (.9,0) -- +(-.15,-.15);
      \draw[dashed] (2,0) circle (.2);
  \end{tikzpicture}
  && $\ovex{A}_1$ \\ \\
  \begin{tikzpicture}[scale=0.6]
    \filldraw[] (0,0) -- +(1.3,0);
      \filldraw[] (2.7,0) -- +(2.3,0);
      \filldraw[] (0,0) -- +(2,1);
        \filldraw[] (4,0) -- +(-2,1);
        \foreach \x in {0,1}
                 {    \filldraw[] (\x,0)  circle (.1cm); }
                 \foreach \x in {3,4,5}
             {    \filldraw[] (\x,0)  circle (.1cm); }
\filldraw[] (2,1) circle (.1cm);
             \draw (2,0) node[] {$\cdots$};
             \draw[dashed] (5,0) circle (.2);
  \end{tikzpicture}
  && $\ovex{A}_n \; (n \geq 2)$\\ \\
  \begin{tikzpicture}[scale=0.6]
    \filldraw[] (0,0) -- +(2.3,0);
      \filldraw[] (3.7,0) -- +(3.3,0);
  \filldraw[] (1,0) -- +(0,0.75);
        \filldraw[] (5,0) -- +(0,0.75);
        \foreach \x in {0,1,2}
                 {    \filldraw[] (\x,0)  circle (.1cm); }
                \foreach \x in {4,5,6,7}
             {    \filldraw[] (\x,0)  circle (.1cm); }
             \filldraw[] (1,0.75) circle (.1cm);
             \filldraw[] (5,0.75) circle (.1cm);
             \draw (3,0) node[] {$\cdots$};
                          \draw[dashed] (7,0) circle (.2);
            \end{tikzpicture} &&
 $ \ovex{D}_{n} \; (n \geq 4)$ \\ \\
  \begin{tikzpicture}[scale=0.6]
  \filldraw[] (0,0) -- +(5,0);
      \filldraw[] (2,0) -- +(0,1.5);
    \foreach \x in {1,...,6} 
             {    \filldraw[] (\x-1,0)  circle (.1cm); }
             \filldraw[] (2,0.75) circle (.1cm);
             \filldraw[] (2,1.5) circle (.1cm);
                          \draw[dashed] (5,0) circle (.2);
  \end{tikzpicture} &&
 $ \ovex{E}_{6}$ \\ \\
  \begin{tikzpicture}[scale=0.6]
  \filldraw[] (-1,0) -- +(7,0);
      \filldraw[] (2,0) -- +(0,1);
    \foreach \x in {1,...,7} 
             {    \filldraw[] (\x-1,0)  circle (.1cm); }
             \filldraw[] (2,1) circle (.1cm);
             \filldraw[] (-1,0) circle (.1cm);
                                      \draw[dashed] (6,0) circle (.2);
              \end{tikzpicture} &&
$\ovex{E}_{7}$ \\ \\
  \begin{tikzpicture}[scale=0.6]
  \filldraw[] (0,0) -- +(8,0);
      \filldraw[] (2,0) -- +(0,1);
    \foreach \x in {1,...,9} 
             {    \filldraw[] (\x-1,0)  circle (.1cm); }
                          \draw[dashed] (8,0) circle (.2);
             \filldraw[] (2,1) circle (.1cm);
  \end{tikzpicture} &&
$\ovex{E}_{8}$
\end{tabular}
\caption{$\Ext$ type diagrams} \label{fig:exts}
\end{figure}

We let $\Ext$ denote the set of overextensions. It is easy to see that the following is the complete list of overextensions, up to isomorphism: \[\ovex{A}_n \;(n \geq 1), \;\; \ovex{D}_n \;(n \geq 4), \; \ovex{E}_n \; (n=6,7,8)\]
(see Figure~\ref{fig:exts}). Here, $\ovex{X}_n$ has $n+2$ vertices. We remark that the corresponding GCMs are all nonsingular; hence a \pisystem of $\Ext$ type is necessarily linearly independent.

\subsection{}\label{sec:finpart-nots} 
From figure~\ref{fig:exts}, one makes the important observation (via case-by-case check) that if $Z$ is an overextension, then the vertex $p$ satisfying the condition in definition \ref{ovex-def} is unique. This vertex is marked by a dashed circle in figure~\ref{fig:exts}. We will call $p$ the {\em overextended vertex} of $Z$, and $Y$ the {\em affine part} of $Z$.

We had $\form{\delta_Y}{\alpha_p} = -1$. Let $\delta_Y = \sum_{q \in Y} c_q \alpha_q$ with $c_q \in \integers_+$ for all $q$.  Observing that $c_q \form{\alpha_q}{\alpha_p} \leq 0$ for all $q$, it follows that: (i)  There is a unique vertex $q$ of $Y$ such that $\form{\alpha_q}{\alpha_p} \neq 0$, (ii) For this vertex, we have $c_q=1$ and $\form{\alpha_q}{\alpha_p} =-1$, (iii) In particular, this means $q$ is a {\em special vertex} of the affine diagram $Y$ (in the terminology of Kac, Chapter 6). Let $\finpart{Z}$ denote the finite type diagram obtained from $Y$ by deleting $q$.
We will call it the {\em finite part} of $Z$. We note that: 
\[ \delta_Y = \alpha_q + \theta_{\finpart{Z}}\]
where $\theta_{\finpart{Z}}$ denotes the highest root of $\finpart{Z}$.
It will be convenient to denote $Y$ by $\affpart{Z}$.

\subsection{}
The following trivial observation is useful: let $X$ be a simply-laced Dynkin diagram and $Z$ a diagram of $\Ext$ type. Suppose there exists $\pi$, a \pisystem of type $Z$ in $X$; we let $\finpart{\pi}, \affpart{\pi}$ denote the subsets of $\pi$ corresponding to the finite and affine parts of $Z$ respectively. For any $w \in W(X)$, $w\pi$ is a \pisystem of type $Z$ in $X$ and $\finpart{(w\pi)} = w (\finpart{\pi}), \; \affpart{w \pi} = w(\affpart{\pi})$.

\subsection{Hyperbolics}
We recall that an indecomposable, {\em symmetrizable} GCM $A$ is said to be of {\em Hyberbolic type} if it is not of finite or affine type and every proper principal submatrix of $A$ is a direct sum of finite or affine type GCMs.

There are  finitely many GCMs of hyperbolic type in ranks 3-10 and infinitely many in rank 2. The former were enumerated, to varying degrees of completeness and detail, in \cite{sacl, buyl-schom,koba-mori}. More recently, this list was organized and independently verified in \cite{lisa-et-al}. We will use this latter reference as our primary source for the Dynkin diagrams of hyperbolic type. Note that \cite{lisa-et-al} does not require symmetrizability in the definition of a hyperbolic type GCM: it contains 142 symmetrizable and 96 non-symmetrizable ones. We let $\Hyp$ denote the set of all symmetrizable GCMs of hyperbolic type of rank $\geq 3$.

    \subsection{} We recall from \S\ref{sec:subdig-po} the {\em subdiagram partial order} on the set of symmetrizable GCMs. We write $B \subseteq A$ if the Dynkin diagram $S(B)$ is a subdiagram of $S(A)$; equivalently $B$ is a principal submatrix of $A$, possibly after a simultaneous permutation of its rows and columns. This is clearly a partial order, once we identify the matrices $\{PAP^T: P \text{ is a permutation matrix}\}$ with each other.

    \subsection{} \label{sec:sl-hyps-list}
\begin{figure}[h]
  \begin{center}
    \begin{tikzpicture}[scale=0.6]
      \filldraw[] (0,0) circle (.1);
      \filldraw[] (1,0) circle (.1);
      \filldraw[] (2,0) circle (.1);
      \filldraw[] (0,-.05) -- +(2,0);
      \filldraw[] (0,.05) -- +(2,0);
      \filldraw[] (.1,0) -- +(.15,.15);
      \filldraw[] (.1,0) -- +(.15,-.15);
      \filldraw[] (.9,0) -- +(-.15,.15);
      \filldraw[] (.9,0) -- +(-.15,-.15);
      \filldraw[] (1.1,0) -- +(.15,.15);
      \filldraw[] (1.1,0) -- +(.15,-.15);
      \filldraw[] (1.9,0) -- +(-.15,.15);
      \filldraw[] (1.9,0) -- +(-.15,-.15);

  \end{tikzpicture}
\hspace{1cm}   \begin{tikzpicture}[scale=0.6]
    \filldraw[] (0,0) circle (.1);
    \filldraw[] (1,0) circle (.1);

     \filldraw[] (0,-.05) -- +(1,0);
      \filldraw[] (0,.05) -- +(1,0);

      \filldraw[] (.1,0) -- +(.2,.15);
      \filldraw[] (.1,0) -- +(.2,-.15);
      \filldraw[] (.9,0) -- +(-.2,.15);
      \filldraw[] (.9,0) -- +(-.2,-.15);

      \begin{scope}[rotate=60]
     \filldraw (1,0) circle (.1);
           \filldraw (0,-.05) -- +(1,0);
      \filldraw (0,.05) -- +(1,0);
     \filldraw (.1,0) -- +(.2,.15);
      \filldraw (.1,0) -- +(.2,-.15);
      \filldraw (.9,0) -- +(-.2,.15);
      \filldraw (.9,0) -- +(-.2,-.15);
      \end{scope}
      
      \begin{scope}[shift={(1,0)}]
        \filldraw[rotate=120] (0,0) --+(1,0);
      \end{scope}
\end{tikzpicture}
\hspace{1cm}
\begin{tikzpicture}[scale=0.6]
  \filldraw[] (0,0) circle (.1);
  \filldraw[] (1,0) circle (.1);
  
  \filldraw[] (0,-.05) -- +(1,0);
  \filldraw[] (0,.05) -- +(1,0);
  
  \filldraw[] (.1,0) -- +(.2,.15);
  \filldraw[] (.1,0) -- +(.2,-.15);
  \filldraw[] (.9,0) -- +(-.2,.15);
  \filldraw[] (.9,0) -- +(-.2,-.15);

  \begin{scope}[rotate=60]
    \filldraw (1,0) circle (.1);
    \filldraw (0,-.05) -- +(1,0);
    \filldraw (0,.05) -- +(1,0);
    \filldraw (.1,0) -- +(.2,.15);
    \filldraw (.1,0) -- +(.2,-.15);
    \filldraw (.9,0) -- +(-.2,.15);
    \filldraw (.9,0) -- +(-.2,-.15);
  \end{scope}

  \begin{scope}[shift={(1,0)},rotate=120]
    \filldraw (1,0) circle (.1);
    \filldraw (0,-.05) -- +(1,0);
    \filldraw (0,.05) -- +(1,0);
    \filldraw (.1,0) -- +(.2,.15);
    \filldraw (.1,0) -- +(.2,-.15);
    \filldraw (.9,0) -- +(-.2,.15);
    \filldraw (.9,0) -- +(-.2,-.15);
  \end{scope}
\end{tikzpicture}
\hspace{1cm}   \begin{tikzpicture}[scale=0.6]
    \filldraw[] (0,0) circle (.1);
    \filldraw[] (1,0) circle (.1);
    \filldraw (0,0) --+(1,0);

      \begin{scope}[rotate=60]
     \filldraw (1,0) circle (.1);
           \filldraw (0,-.05) -- +(1,0);
      \filldraw (0,.05) -- +(1,0);
     \filldraw (.1,0) -- +(.2,.15);
      \filldraw (.1,0) -- +(.2,-.15);
      \filldraw (.9,0) -- +(-.2,.15);
      \filldraw (.9,0) -- +(-.2,-.15);
      \end{scope}
      
      \begin{scope}[shift={(1,0)}]
        \filldraw[rotate=120] (0,0) --+(1,0);
      \end{scope}
\end{tikzpicture}

\vspace{1cm}
\begin{tikzpicture}[scale=0.5]

\foreach \x in {0,1,2}
{\begin{scope}[rotate=120*\x]
\begin{scope}[rotate=0]
\filldraw[] (0,0) circle (.1);
    \filldraw[] (1,0) circle (.1);
    \filldraw (0,0) --+(1,0);
        \coordinate (A) at (1,0);
\end{scope}

\begin{scope}[rotate=120]
\filldraw[] (0,0) circle (.1);
    \filldraw[] (1,0) circle (.1);
    \filldraw (0,0) --+(1,0);
    \coordinate (B) at (1,0);
\end{scope}

\draw (A) -- (B);
\end{scope}}
\end{tikzpicture}
\hspace{1cm}
\begin{tikzpicture}[scale=0.5]

\foreach \x in {0,1}
{\begin{scope}[rotate=120*\x]
\begin{scope}[rotate=0]
\filldraw[] (0,0) circle (.1);
    \filldraw[] (1,0) circle (.1);
    \filldraw (0,0) --+(1,0);
        \coordinate (A) at (1,0);
\end{scope}

\begin{scope}[rotate=120]
\filldraw[] (0,0) circle (.1);
    \filldraw[] (1,0) circle (.1);
    \filldraw (0,0) --+(1,0);
    \coordinate (B) at (1,0);
\end{scope}

\draw (A) -- (B);
\end{scope}}
\end{tikzpicture}
\hspace{1cm}
\begin{tikzpicture}[scale=0.4]
  \filldraw (0,0) circle (.1cm);
  \filldraw (2,0) circle (.1cm);
  \filldraw (0,2) circle (.1cm);
  \filldraw (2,2) circle (.1cm);
  \filldraw (1,1) circle (.1cm);
  
  \draw (0,0) -- (2,0) -- (2,2) -- (0,2) -- cycle;
  \draw (2,0) -- (0,2);
  
\end{tikzpicture}
\hspace{1cm}
\begin{tikzpicture}[scale=0.5]
  \foreach \x in {0,...,4}
{  \begin{scope}[rotate=72*\x]
    \filldraw (0,0) circle (.1cm);
    \filldraw (1,0) circle (.1cm);
    \filldraw (0,0) -- +(1,0);
\end{scope}
}
\end{tikzpicture}
\end{center}
\caption{Simply-laced hyperbolics (ranks 3-10) that are not of $\Ext$ type.}\label{fig:irreg-hyps}
\end{figure}
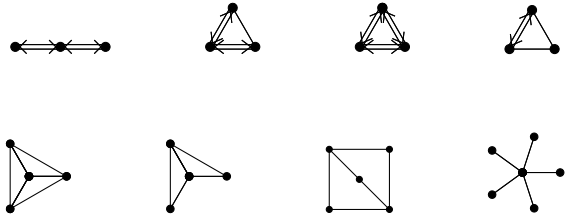

We now isolate the {\em symmetric} GCMs of hyperbolic type. By checking the classification case-by-case (see for instance \cite[Tables 1,2]{svis-e10} or \cite{lisa-et-al}), one finds that these are either (i) of $\Ext$ type:
\begin{equation} \label{eq:hyp-ext}
  \ovex{A}_n, \; (1 \leq n \leq 7), \;\; \ovex{D}_n, \; (4 \leq n \leq 8), \;\; \ovex{E}_n, \; (6 \leq n \leq 8)
\end{equation}
or (ii) one of the diagrams in Figure~\ref{fig:irreg-hyps}, or (iii) one of the rank 2 symmetric GCMs $\begin{bmatrix} 2 & -a \\ -a & 2 \end{bmatrix}$ for $a \geq 3$. We observe by inspection of figure~\ref{fig:exts} that the diagrams in (ii) and (iii) do not contain a subdiagram of $\Ext$ type.

\subsection{}
 The next lemma underscores the special role played by the hyperbolic overextensions. These are precisely the minimal elements of the set of overextensions relative to the partial order $\subseteq$.
\begin{lemma} \label{lem:minelts-ext}
$$\min (\Ext, \subseteq) = \Ext \cap \Hyp$$
\end{lemma}
\begin{proof}
Observe that $\ovex{E}_7 \subseteq \ovex{A}_n$ for $n \geq 8$ and $\ovex{E}_8 \subseteq \ovex{D}_n$ for $n \geq 9$. We are thus left with the diagrams of equation~\eqref{eq:hyp-ext} as possible candidates for minimal elements. Now, each of these diagrams except $\ovex{D}_8$ contains a unique subdiagram of affine type, obtained by removing a single vertex. So these diagrams cannot contain a proper subdiagram of $\Ext$ type. As for the diagram $Z=\ovex{D}_8$, it contains two subdiagrams of affine type, $Y_1 = E_8^{(1)}$ and $Y_2 = D_8^{(1)}$, obtained by deleting appropriate vertices $p_1, p_2$,  but only the former satisfies $\form{\delta_Y}{\alpha_p}=-1$ (this is $-2$ for the latter). Thus, $\ovex{D}_8$ is also minimal. 
\end{proof}

\section{Weyl group orbits of \pisystems of type $\ovex{A}_1$}
In this section, we focus on the diagram $\ovex{A}_1$. The corresponding Kac-Moody algebra was first studied by Feingold and Frenkel \cite{feingold-frenkel}.
\subsection{}
We consider the problem of determining $\mult(\ovex{A}_1, X)$ for a simply-laced Dynkin diagram $X$. This is an important special case of the more general result of the next section. The latter result will be obtained by arguments similar to the ones used here, albeit with more notational complexity.

\subsection{}
We begin with the following lemma which asserts that every Dynkin diagram of $\Ext$ type has a ``canonical''  \pisystem of type $\hasl$.

\begin{lemma}
  Given a Dynkin diagram $Z$ of $\Ext$ type, define:
  \[\pi(Z):=\{\theta_{\finpart{Z}}, \delta_Y - \theta_{\finpart{Z}}, \alpha_p\}\]
  (notations $\finpart{Z}, Y, \,p,  \,\theta_{\finpart{Z}}$ are as defined in \S\ref{sec:finpart-nots}). Then $\pi(Z)$ is a linearly independent, positive  \pisystem of type $\hasl$.
\end{lemma}

\begin{proof}
  We only need to show that the type of $\pi(Z)$ is $\hasl$, the other assertions following from the observation that the three roots in $\pi(Z)$ are real, positive and have disjoint supports (cf. \S\ref{sec:finpart-nots}). Since $Z$ is simply-laced, we normalize the form such that all real roots have norm 2. Thus
  \[\form{\theta_{\finpart{Z}}}{\delta_Y - \theta_{\finpart{Z}}} = - \form{\theta_{\finpart{Z}}}{\theta_{\finpart{Z}}}  = -2\]
  It is clear from \S\ref{sec:finpart-nots} that $\form{\theta_{\finpart{Z}}}{\alpha_p} =0$ and $\form{\delta_Y}{\alpha_p} =-1$. This completes the verification.
\end{proof}

\bthm \label{thm:ha11-pi-systems}
Let $X$ be a simply-laced Dynkin diagram. Then:
\be
\item $X$ has a \pisystem of type $\hasl$ if and only if it contains a subdiagram of $\Ext$ type.
\item The number of $W(X)$-orbits of \pisystems of type $\hasl$ in $X$ is twice the number of such subdiagrams (and is, in particular, finite).
\ee
\ethm

\begin{proof}
  In light of Theorem~\ref{thm:pos-neg}, any \pisystem of type $\hasl$ in $X$ is $W(X)$-equivalent to a positive or a negative \pisystem, but not both. Thus, to prove the above theorem, it is sufficient to construct a bijection from the set of $\Ext$ type subdiagrams of $X$ to $W(X)$-equivalence classes of {\em positive} \pisystems  of type $\hasl$ in $X$. We claim that the following map defines such a bijection: \[Z \mapsto [\pi(Z)]\] 

  We will first establish the injectivity. Suppose $Z_1, Z_2$ are $\Ext$ type subdiagrams of $X$, with affine parts $Y_1, Y_2$ and overextended vertices $p_1, p_2$ respectively. Suppose $\pi(Z_1) \sim \pi(Z_2)$ i.e., there exists $\sigma \in W(X)$ such that $\sigma (\pi(Z_1)) = \pi(Z_2)$.
  Consider the \pisystems: \[\pi_j = \{\theta_{\finpart{Z}_j}, \delta_{Y_j} - \theta_{\finpart{Z}_j}\}, \;\;\;j=1, 2.\]
  We note that:
  \be
  \item $\pi_j$ is of type $\afsl$.
  \item $\pi_j$ is supported in the affine subdiagram $Y_j$ of $X$.
  \item $\sigma(\pi_1) = \pi_2$.
  \ee
  Now, it follows from part (2) of theorem~\ref{thm:affine-pi} that $Y_1 = Y_2$ and $\sigma \in W(Y_1 \sqcup Y_1^\perp)$. Since $p_1 \not\in Y_1 \sqcup Y_1^\perp$, we can only have $\sigma \alpha_{p_1} = \alpha_{p_2}$ if $p_1 = p_2$. Thus, $Z_1 = Z_2$ as required.

  \medskip
  Next, we turn to the surjectivity of this map.
  Let $\{\beta_{-1}, \beta_0, \beta_1\}$ be a positive \pisystem of $X$ of type $\hasl$. Since $\{\beta_0, \beta_1\}$ form a \pisystem of type $A_1^{(1)}$, which is affine, it follows from theorem~\ref{thm:affine-pi} that there is  a unique affine type subdiagram $Y$ of $X$ and an element $w \in W(X)$ such that $w\beta_i$ is supported in $Y$ for $i=0,1$. Further (as in the proof of theorem~\ref{thm:affine-pi}), since $w(\beta_0 + \beta_1)$ is an isotropic root of $\lieg(X)$, we must have
  $w(\beta_0 + \beta_1) = k\delta_Y$ for some {\em nonzero} integer $k$. Since $\form{\beta_0 + \beta_1}{\beta_{-1}} = -1$, we conclude $k=\pm 1$.
But $\beta_0 + \beta_1 \in Q_+(X)$ by proposition \ref{prop:criterion-pos-neg}, and $w^{-1}(\delta_Y) \in \imroots[+]$ since $\delta_Y$ is a positive imaginary root. This implies $k=1$.

\medskip
Let $\beta'_i = w\beta_i$; thus $\beta'_0, \beta'_1$ are supported in $Y$, their sum equals $\delta_Y$ and $\form{\delta_Y}{\beta'_{-1}} = -1$.
We now need the following lemma:
\begin{lemma} \label{lem:oshima}
  Let $X$ be a simply-laced Dynkin diagram, $Y$ an affine subdiagram of $X$ and $\beta$ a real root of $X$ satisfying $\form{\delta_Y}{\beta} = -1$. Then there exists $\sigma \in W(Y \sqcup Y^\perp)$ such that $\sigma \beta$ is a simple root of $X$.
\end{lemma}

We defer the proof of this lemma to the next subsection. Here, we use it to complete the proof of Theorem~\ref{thm:ha11-pi-systems}. We take $\beta = \beta'_{-1}$ in lemma~\ref{lem:oshima}. We obtain $\sigma \in W(Y \sqcup Y^\perp)$ such that $\sigma \beta'_{-1} = \alpha_p$ for some vertex $p$ of $X$.  Define $Z := Y \cup \{p\}$. Since $\sigma$ stabilizes $\delta_Y$, we have $\form{\delta_Y}{\alpha_p}=-1$; thus $Z$ is of $\Ext$ type.

    Since $\beta'_0, \beta'_1$ are supported in $Y$, so are $\sigma \beta'_0,  \sigma \beta'_1$; further $\sigma \beta'_0 + \sigma \beta'_1 = \delta_Y$. Now
    \[\form{\sigma \beta'_1}{\alpha_p} = \form{\sigma \beta'_1}{\sigma \beta'_{-1}} = 0.\]
   This implies that $\sigma \beta'_1$ is supported in $\finpart{Z}$.
    Since $\finpart{Z}$ is a simply-laced finite type diagram, all its real roots are conjugate under its Weyl group.
    Thus, there exists $\tau \in W(\finpart{Z})$ such that $\tau\sigma \beta'_1= \theta_{\finpart{Z}}$.
Since $\tau$ stabilizes both $\delta_Y$ and $\alpha_p$, we conclude that $\{\tau\sigma \beta'_i: i=-1,0,1\} = \pi(Z)$, as required. 
\end{proof}

\subsection{} We now turn to the proof of Lemma~\ref{lem:oshima}. We use the notations of the lemma.
Since $\delta_Y$ is an antidominant weight of $X$, $\beta$ must be a positive root. Further it is clear from $\form{\delta_Y}{\beta}=-1$ that $\beta$ 
   must have the form:
   \begin{equation} \label{eq:beta-form}
     \beta = \alpha_p + \sum_{q \in Y \sqcup Y^\perp} c_q(\beta) \alpha_q
   \end{equation}
where $p$ is a vertex of $X$ such that $\form{\delta_Y}{\alpha_p} =-1$, and $c_q(\beta)$ are non-negative integers. Consider the $W(Y \sqcup Y^\perp)$-orbit of $\beta$. Since the coefficient of $\alpha_p$ remains the same, any element $\gamma$ of this orbit is a positive root that has the same form as the right hand side of \eqref{eq:beta-form} for some non-negative coefficients $c_q(\gamma)$. Let $\gamma$ be a minimal height element of this orbit, i.e., one for which $\sum_q c_q(\gamma)$ is minimal. Then, we have: (i) $\form{\gamma}{\alpha_q} \leq 0$ for all $q \in Y \sqcup Y^\perp$, since otherwise $s_q \gamma$ would have strictly smaller height, (ii) $\form{\gamma}{\gamma} = \form{\alpha_p}{\alpha_p}$ since all real roots have the same norm ($X$ is simply-laced). We compute:
\[ 0 = \form{\gamma + \alpha_p}{\gamma - \alpha_p} = \sum_{q \in Y \sqcup Y^\perp} c_q(\gamma) \, \form{\gamma + \alpha_p}{\alpha_q}\]
Since $\form{\alpha_p}{\alpha_q} \leq 0$, we conclude from (i) above that either $ c_q(\gamma) =0$ or $\form{\gamma}{\alpha_q} =  \form{\alpha_p}{\alpha_q} =0$ for each $q \in  Y \sqcup Y^\perp$. If some $c_q(\gamma) \neq 0$, it would imply that $\gamma$ has disconnected support, which is impossible since $\gamma$ is a root. Thus, $\gamma = \alpha_p$ and the proof of the lemma is complete. \qed

\subsection{}
We note that the key step in the proof above was showing that the set of all real roots $\beta$ which have the form of equation~\eqref{eq:beta-form} forms a single orbit under the standard parabolic subgroup $W(Y \sqcup Y^\perp)$ of $W$. In fact, those very same arguments prove a strengthened assertion. We formulate this below.

Given a Dynkin diagram $X$ with simple roots $\alpha_i$ and given any $\alpha$ in its root lattice, we define the coefficients $c_i(\alpha)$ by:
\[ \alpha = \sum_{i \in X} c_i(\alpha) \, \alpha_i \]
If $J$ is a subdiagram of $X$, we define $\alpha_J = \sum_{i \in J} c_i(\alpha)\,\alpha_i$ and 
$\alpha_J^\dag = \sum_{i \not\in J} c_i(\alpha)\,\alpha_i$.

\begin{proposition}\label{prop:oshima-gen}
  Let $X$ be a symmetrizable Dynkin diagram with invariant bilinear form $\form{\cdot}{\cdot}$ and simple roots $\alpha_i$. Let $J$ be a subdiagram of $X$, and fix a {\em nonzero} element $\zeta = \sum_{i \not\in J} b_i\, \alpha_i$ of the root lattice of $X\backslash J$.  Consider the set
  \[\mathcal{O} = \{\beta \in \reroots(X): \beta_J^\dag = \zeta \text{ and } \form{\beta}{\beta} = \form{\zeta}{\zeta}\}\]
  Then:
  \be
  \item If $\zeta$ is a root of $\lieg(X\backslash J)$, then $\mathcal{O} = W_J \,\zeta$ where $W_J$ is the standard parabolic subgroup
    $\langle s_j: j \in J\rangle$ of $W$.
  \item If $\zeta$ is not a root of $\lieg(X \backslash J)$, then $\mathcal{O}$ is empty.
  \ee
\end{proposition}

\begin{proof}
Suppose $\mathcal{O}$ is non-empty, then $\zeta$ or $-\zeta$ lies in $Q_+(X\backslash J)$. We may assume the former case holds, so in fact
$\mathcal{O} \subset \reroots[+](X)$. Since $\mathcal{O}$ is $W_J$-stable, it decomposes into $W_J$-orbits. Let $\mathcal{O}'$ denote one such orbit. let $\beta$ denote an element of minimal height in $\mathcal{O}'$; as in the proof of Lemma~\ref{lem:oshima}, this implies $\form{\beta}{\alpha_j} \leq 0$ for all $j \in J$; hence $\form{\beta}{\alpha} \leq 0$ for all elements $\alpha \in Q_+(J)$.
We now have $0 = \form{\beta+\zeta}{\beta-\zeta} = \form{\beta+\beta_J^\dag}{\beta_J}$. But as observed already, $\form{\beta}{\beta_J} \leq 0$; further $\form{\beta_J^\dag}{\beta_J} \leq 0$ since these elements have disjoint supports. This implies $\form{\beta}{\beta_J} = \form{\beta_J^\dag}{\beta_J} =0$. Suppose $\beta_J$ is nonzero, the latter implies that $\beta = \beta_J + \beta_J^\dag$ has disconnected support. Hence it cannot be a root. This contradiction shows $\beta_J=0$, i.e., $\beta = \beta_J^\dag = \zeta$. In particular, $\zeta$ is a root, and belongs to any $W_J$ orbit in $\mathcal{O}$. Hence $\mathcal{O} = W_J\,\zeta$.
\end{proof}

\bremark
\be
\item If $X$ is simply-laced and $J$ is a singleton, say $J=\{p\}$, and $\zeta = \alpha_p$, then $\mathcal{O}$  consists precisely of those real roots $\beta$ of $X$ which have the form of equation~\eqref{eq:beta-form}.
\item If $X$ is of finite type and $\zeta$ is a root of $X\backslash J$, then Proposition~\ref{prop:oshima-gen} is a consequence of {\em Oshima's lemma} \cite[Lemma 4.3]{oshima}, \cite[Lemma 1.2]{dyer-lehrer}.
\ee
\eremark

\subsection{}
We now have the following corollary of Theorem~\ref{thm:ha11-pi-systems}.

\bcor \label{cor:hyp-a11}
Let $X$ be a Dynkin diagram of $\Ext$ type. Then:
\be
\item If $X \in \Hyp$, then there are exactly two \pisystems of type $\hasl$ in $X$, up to $W(X)$-equivalence. In other words:
  \[\mult(\hasl,X)=2 \text{ for } X=\ovex{A}_n \, (1 \leq n \leq 7), \; \ovex{D}_n \, (4 \leq n \leq 8), \; \ovex{E}_n \, (n=6,7,8).\]
\item $\mult(\ovex{A}_1, \ovex{A}_8) = 6$, $\mult(\ovex{A}_1, \ovex{A}_n) = 10$ for $n \geq 9$.
\medskip
\item $\mult(\ovex{A}_1, \ovex{D}_9) = 6$, $\mult(\ovex{A}_1, \ovex{D}_n) = 4$ for $n \geq 10$.
\ee
\ecor
\bmyproof
The first part follows from Lemma~\ref{lem:minelts-ext} and Theorem~\ref{thm:ha11-pi-systems}.
For parts (2), (3), we need to count the number of subdiagrams of the ambient diagram which are of $\Ext$ type. We list these out in each case, leaving the easy verification to the reader.
\be
\item  $\ovex{A}_8$: one subdiagram of type $\ovex{A}_8$ and two of type $\ovex{E}_7$.
\smallskip\item  $\ovex{A}_n \; (n \geq 9)$: one subdiagram of type $\ovex{A}_n$ and two each of types $\ovex{E}_7$ and $\ovex{E}_8$.
\smallskip\item  $\ovex{D}_9$:  one subdiagram of type $\ovex{D}_9$ and two of type $\ovex{E}_8$.
\smallskip\item  $\ovex{D}_n \; (n \geq 10)$:  one subdiagram of type $\ovex{D}_n$ and one of type $\ovex{E}_8$.
\ee
\emyproof

We also have the following result concerning the simply-laced hyperbolic diagrams not included in the previous corollary.
\bcor
Let $X$ be a simply-laced hyperbolic Dynkin diagram. If $X \not\in \Ext$, then $X$ does not contain a \pisystem of type $\hasl$.
\ecor
\bmyproof
This follows from the observation made in \S\ref{sec:sl-hyps-list} that such diagrams do not contain subdiagrams of $\Ext$ type.
\emyproof

Finally, we remark that Theorem~\ref{thm:ha11-pi-systems} can be applied just as easily even when $X$ is neither in $\Ext$ nor $\Hyp$. For example, the diagram $X = E_{11}$, obtained by further extension of $\ovex{E}_8$ \cite{Henneaux2008} contains a unique subdiagram of $\Ext$ type, namely $\ovex{E}_8$. Thus, $\mult(\ovex{A}_1, E_{11}) = 2$.

\section{The general case}

\bthm \label{thm:mainthm-gen}
Let $X$ be a simply-laced Dynkin diagram and let $K$ be a diagram  of $\Ext$ type. Then:

\be
\item There exists a \pisystem in $X$ of type $K$ if and only if there exists an $\Ext$ type subdiagram $Z$ of $X$ such that $\finpart{Z}$ has a \pisystem of type $\finpart{K}$. 

\item The number of $W(X)$ orbits of \pisystems of type $K$ in $X$ is given by: 

  \begin{equation} \label{eq:gen-case-mults}
    \mult{(K,X)} =  2 \sum_{\substack{Z \subseteq X \\ Z \in \Ext}} \mult(\finpart{K}, \finpart{Z})
  \end{equation}
\ee

where $\finpart{K}, \finpart{Z}$ denote their finite parts.
\ethm

We remark that equation \eqref{eq:gen-case-mults} reduces the computation of the multiplicity of $K$ in $X$ to a sum of multiplicities involving only finite type diagrams. The latter, as mentioned earlier, are completely known \cite{dynkin}. Observe also that for $K = \ovex{A}_1$, $\finpart{K}$ is of type $A_1$. Since any $\finpart{Z}$ occurring on the right hand side of \eqref{eq:gen-case-mults} is simply-laced, we have $\mult(\finpart{K}, \finpart{Z}) =1$. So this reduces exactly to Theorem \ref{thm:ha11-pi-systems} in this case.

\begin{corollary}\label{cor:gencase-cor} Let $K$ be a Dynkin diagram of $\Ext$ type. Then,
\be
\item $\mult(K,X)$ is finite for all simply-laced diagrams $X$.
\item $\mult(K,X) = 2 \mult(\finpart{K},\finpart{X})$ for all $X \in \Hyp \cap \Ext$. 
\ee
\end{corollary}

We now prove theorem \ref{thm:mainthm-gen}.

\begin{proof}
It is enough to prove the second part of the theorem. Now, by Theorem~\ref{thm:pos-neg}, any \pisystem in $X$ of type $K$ is either positive or negative, but not both. Consider the sets:
\begin{itemize}
\setlength\itemsep{1em}
\item $\mathcal{A}$: the set of $W(X)$-orbits of positive \pisystems of type $K$ in $X$;
\item $\widehat{\mathcal{B}}$: the set of all pairs $(Z,\Sigma)$ where $Z$ is an $\Ext$ type subdiagram of $X$ and $\Sigma$ is a  positive
 \pisystem of type $\finpart{K}$ in $\finpart{Z}$.
\item $\mathcal{B} =  \widehat{\mathcal{B}}/\!\!\sim$, the equivalence classes of $\widehat{\mathcal{B}}$ under the equivalence relation defined by: $(Z, \Sigma) \sim (Z', \Sigma')$ if and only if $Z=Z'$ and $\Sigma'$ is in the $W(\finpart{Z})$-orbit of $\Sigma$.
\end{itemize}

\medskip
Since $2|\mathcal{A}|$ and $2|\mathcal{B}|$ are the two sides of equation \eqref{eq:gen-case-mults}, it is sufficient to construct a bijection from the set $\mathcal{B}$ to $\mathcal{A}$. We first define a map from $\widehat{\mathcal{B}}$ to $\mathcal{A}$. Let
$(Z, \Sigma) \in \widehat{\mathcal{B}}$. Let $\finpart{Z}$ and $\affpart{Z}$ denote the finite and affine parts of $Z$, and let $p$ denote its  overextended vertex. Since $\Sigma$ is a \pisystem of type $\finpart{K}$ in $\finpart{Z}$, we identify $\Delta(\finpart{K})$ with a subset of $\Delta(\finpart{Z})$ via corollary~\ref{cor:real-to-real}. Let $\theta_\Sigma$ denote the highest root in $\Delta(\finpart{K})$ (identified with its image in $\Delta(\finpart{Z}) \subset Q(Z)$). Consider the set
\[\pi(Z,\Sigma)=\{\alpha_{p},\delta_{\affpart{Z}}-\theta_{\Sigma}\} \cup \Sigma \]
It is straighforward to see that this is a \pisystem. Further, it is of type $K$.
We now claim that the map: $\widehat{\mathcal{B}} \to \mathcal{A}, \;\;(Z,\Sigma) \mapsto [\pi(Z,\Sigma)]$ factors through $\mathcal{B}$ and defines a bijection between $\mathcal{B}$ and $\mathcal{A}$.

Firstly, suppose $(Z,\Sigma) \sim (Z,\Sigma')$, i.e., $w\Sigma = \Sigma'$ for some $w \in W(\finpart{Z})$. Since clearly $w\alpha_p = \alpha_p$,
$w\delta_{\affpart{Z}} = \delta_{\affpart{Z}}$ and $\theta_{\Sigma'} = w\theta_{\Sigma}$, we conclude that $\pi(Z,\Sigma') = w\,\pi(Z,\Sigma)$. So the map does indeed factor through $\mathcal{B}$. We will now show it is an injection.

Suppose $(Z_i,\Sigma_i) \in \widehat{\mathcal{B}}$, $i=1,2$ are such that $[\pi(Z_1,\Sigma_1)] = [\pi(Z_2,\Sigma_2)]$, i.e.,  there exists $\sigma \in W(X)$ such that $\sigma (\pi(Z_1,\Sigma_1)) = \pi(Z_2,\Sigma_2)$. Let $p_i$ denote the overextended vertex of $Z_i$. 
  Consider the \pisystems: \[\pi_j = \{\delta_{\affpart{Z_j}}-\theta_{\Sigma_j}\} \cup \Sigma_j, \;\;\;j=1, 2.\]
  We note that:
  (i) $\pi_j$ is of type $\affpart{K}$,
  (ii) $\pi_j$ is supported in the affine subdiagram $\affpart{Z_j}$ of $X$, and
  (iii) $\sigma(\pi_1) = \pi_2$.

  Now, it follows from part (2) of theorem~\ref{thm:affine-pi} that $\affpart{Z_1} =\affpart{Z_2}$ and $\sigma \in W(\affpart{Z_1} \sqcup \affpart{Z_1}^\perp)$. Since $p_1 \not\in \affpart{Z_1} \sqcup \affpart{Z_1}^\perp$, we can only have $\sigma \alpha_{p_1} = \alpha_{p_2}$ if $p_1 = p_2$. Thus, $Z_1 = Z_2$. We write $\sigma=\tau\tau'$ with $\tau\in W(\affpart{Z_1})$ and $\tau'\in W(\affpart{Z_1}^\perp$).
Since $\sigma \pi_1 = \pi_2$, we obtain $\tau\,\Sigma_1=\Sigma_2$ (in fact, $\tau\,\pi_1=\pi_2$) since $\tau'$ fixes each element of $\pi_1$ pointwise. 
Further, $\sigma \alpha_{p_1} = \alpha_{p_1}$ implies that $\sigma \in W(\{p_1\}^\perp)$. In particular, $\tau \in W(\affpart{Z_1}) \cap W(\{p_1\}^\perp) = W(\finpart{Z_1})$. Hence we obtain $(Z_1,\Sigma_1) \sim (Z_2,\Sigma_2)$, in other words, the map defined above is injective on $\mathcal{B}$.
Next, we show surjectivity of the map. Let $\pi$  be a positive \pisystem in $X$ of type $K$; we will show that $[\pi]$ is in the image of the map. Let $\finpart{\pi}, \affpart{\pi}$ be the subsets of $\pi$ corresponding to the finite and affine parts of $K$ respectively. Now, $\affpart{\pi}$ is a positive \pisystem of type $\affpart{K}$ in $X$. By theorem~\ref{thm:affine-pi}, there is an affine type subdiagram $Y$ of $X$, and an element $w \in W(X)$ such that every element of (the positive \pisystem) $w(\affpart{\pi}) = \affpart{(w\pi)}$ is supported in $Y$. Since $[\pi] = [w\pi]$, let us replace $\pi$ with $w\pi$ in what follows. Thus, $\pi$  is a positive \pisystem of type $K$ such that $\affpart{\pi}$ is supported in $Y$. Let $\beta \in \pi$ correspond to the overextended vertex of $K$, and let $\delta_{\affpart{\pi}}$ denote the null root of $\affpart{K}$, identified with its image in $\Delta({\affpart{\pi}}) \subset \Delta(X)$. Thus $\delta_{\affpart{\pi}}$ (i) is a positive imaginary root of $X$ (by corollary~\ref{cor:real-to-real}), (ii) is supported in $Y$, and (iii) satisfies $\form{\delta_{\affpart{\pi}}}{\beta} = -1$.
The first two conditions imply $\delta_{\affpart{\pi}} = r \delta_{Y}$ for some $r \geq 1$, while the third implies $r=1$.

As in the proof of Theorem \ref{thm:ha11-pi-systems}, we now appeal to Lemma \ref{lem:oshima} to find an element $\sigma \in W(Y \sqcup Y^\perp)$ such that $\sigma \beta = \alpha_p$ for some vertex $p$ of $X$. Define $Z = Y \cup \{p\}$; this is clearly an $\Ext$ type subdiagram of $X$. Consider the positive \pisystem $\xi = \sigma \pi$ of type $K$. We have: (a) $\alpha_p \in \xi$, (b) $\affpart{\xi}$ is supported in $Y$ and (c) $\delta_{\affpart{\xi}} = \delta_Y$. 
Further, $\form{\alpha}{\beta} = 0$ for all $\alpha \in \finpart{\pi}$ gives us $\form{\sigma \alpha}{\sigma \beta} =0 $, i.e., $\form{\alpha'}{\alpha_p} =0$ for all $\alpha' \in \finpart{\xi}$. This in turn implies that: (d) $\finpart{\xi}$ is  supported in $\finpart{Z}$.
From (a), (c) and (d) we conclude $\xi = \pi(Z,\finpart{\xi})$. Since $[\pi] = [\xi]$ and $\finpart{\xi}$ is of type $\finpart{K}$, the proof is complete.

\end{proof}

%% file: physics-embeddings.tex
We make some remarks on the motivations from physics that led to this work. 
It is known that  $E_{10}$ symmetry appears in 11 dimensional supergravity in several ways. In dimensionally reduced supergravity theories, there is a Lie group $G$ and subgroup $K$ such that the coset
space $G/K$ fibers over spacetime $\mathcal{M}$ and the scalar fields  are maps $\Sigma: \mathcal{M}\longrightarrow G/K$.  Under dimensional reduction to dimension $D=1$, $G=E_{10}(\mathbb{R})$ and $K$ is the subgroup fixed by the Cartan involution \cite{dhn,julia}. The real roots and the representations of  $E_{10}$ have been shown to correspond to the fields of 11 dimensional supergravity at low levels \cite{dhn}. 

There is a similar description of the symmetries of Einstein gravity in $D = 4$ spacetime dimensions in terms of the rank 3 Feingold-Frenkel hyperbolic Kac-Moody algebra $A_1^{++}$, also denoted $AE_3$ in the physics literature. Under dimensional reduction to $D=1$ spacetime dimensions, $AE_3$ is conjectured to be a symmetry of the dimensionally reduced Lagrangian \cite{julia}. The real roots and the representations of  $AE_3$  correspond to the fields of gravity at low levels
\cite{dhn}. 

The Lie algebra $E_{10}$ contains a regular subalgebra isomorphic to $AE_3$ (by the results of \cite{svis-e10} or Theorem~\ref{thm:ha11-pi-systems}). This reflects the inclusion of Einstein gravity into $11$-dimensional supergravity. Theorem~\ref{thm:ha11-pi-systems} states that all \pisystems of type $AE_3$ in $E_{10}$  are conjugate (up to negation) under the Weyl group of $E_{10}$, indicating that there is a `canonical' inclusion of Einstein gravity into $11$-dimensional supergravity.  
We remark that there is also a physics-suggested way to embed $AE_3$ into $E_{10}$ via the `gravity truncation' method of \cite[\S 4]{DHi}.



Weyl orbits of \pisystems of affine type occur in \cite[\S 3]{EHKNT}, where the authors describe a family of embeddings $A_1^{(1)} \preceq E_{8}^{(1)} \; (\subseteq E_{10})$ and give a brane interpretation of these.